\documentclass{iopjournal}
\usepackage{amsmath}
\usepackage{bm}
\usepackage{float}
\begin{document}

\articletype{Article type} 

\title{Multigrid Poisson Solver for Complex Geometries Using Finite Difference
Method}

\author{Deepak Gautam$^1$\orcid{0000-0000-0000-0000}, Bhooshan Paradkar$^1$\orcid{0000-0000-0000-0000}}

\affil{$^1$School of Physical Sciences, UM-DAE Centre for Excellence in Basic Sciences, University of Mumbai - 400098,
India}

\affil{$^*$Author to whom any correspondence should be addressed. Email: deepak.gautam@cbs.ac.in}

\email{ deepak.gautam@cbs.ac.in, bhooshan.paradkar@cbs.ac.in}

\keywords{sample term, sample term, sample term}

\begin{abstract}
\vspace{5mm}
We present an efficient numerical method, inspired by transformation optics, for solving the Poisson equation in complex and arbitrarily shaped geometries. The approach operates by mapping the physical domain to a uniform computational domain through coordinate transformations, which can be applied either to the entire domain or selectively to specific boundaries inside the domain. This flexibility allows both homogeneous (Laplace equation) and inhomogeneous (Poisson equation) problems to be solved efficiently using iterative or fast direct solvers, with only the material parameters and source terms modified according to the transformation. The method is formulated within a finite difference framework, where the modified material properties are computed from the coordinate transformation equations, either analytically or numerically. This enables accurate treatment of arbitrary geometric shapes while retaining the simplicity of a uniform grid solver. Numerical experiments confirm that the method achieves second-order accuracy, $\mathcal{O}(\Delta x^{2})$, and offers a straightforward pathway to integrate fast solvers such as multigrid methods on the uniform computational grid.
\end{abstract}

\section{Introduction}
\vspace{5mm}
An efficient and accurate solution of the Poisson equation is fundamental to a wide range of problems in computational physics ~\cite{briggs2000multigrid,trottenberg2001multigrid,hockney1988computer}. In computational plasma physics, the Poisson solvers are particularly important for both kinetic and fluid modeling of plasma. Such solvers are at the heart of electrostatic particle-in-cell (PIC) codes, used for kinetic modeling of plasma. Similarly, the pressure Poisson equation is used in incompressible fluid modeling to calculate the pressure field while maintaining the divergence free velocity field. In large-scale realistic simulations, the solution domain often exhibits non-uniform geometries or mesh distributions, either due to complex physical boundaries or adaptive mesh refinement strategies. For such cases, high-performance solvers must not only maintain accuracy but also remain computationally efficient.\\
Multigrid methods have long been recognized as one of the most effective approaches for solving elliptic equations, offering optimal or near-optimal complexity for a wide range of discretizations~\cite{brandt1977multigrid,hackbusch1985multigrid,xu1992iterative}. In particular, geometric multigrid method can deliver rapid convergence by exploiting the hierarchical structure of the underlying mesh. However, when the computational grid is non-uniform, direct application of geometric multigrid becomes challenging. Finite element methods (FEM) are often employed in such situations because they naturally accommodate unstructured meshes~\cite{zienkiewicz2005fem,elman2014fem}, but their implementation can be complex, and the resulting matrix operations tend to occurs higher computational cost compared to structured finite difference schemes.\\
In contrast, finite difference methods are simple to implement and highly efficient on uniform meshes, making them attractive for large-scale problems with structured geometry. Unfortunately, for non-uniform physical domains, applying finite difference multigrid directly is not straightforward without loss of accuracy or convergence efficiency. An alternative strategy, inspired by transformation optics~\cite{pendry2006science,leonhardt2006science,li2008prl}, is to map the non-uniform physical domain into a uniform computational domain through a suitable coordinate transformation. In transformation optics, Maxwell’s equations in a curved or graded medium are equivalently expressed in a uniform space with spatially varying material properties (permittivity and permeability). These ideas are widely used to design novel materials, typically called meta-materials. By analogy, for the Poisson equation, a similar transformation modifies the permittivity distribution in the computational domain to preserve the original physics.\\
In this work, we exploit this concept to develop a finite difference geometric multigrid solver that operates entirely on a uniform computational grid, even when the underlying physical domain is non-uniform. The transformation from the physical coordinates to the computational coordinates is used to derive an equivalent Poisson equation with a spatially varying permittivity tensor. This approach eliminates the need for FEM assembly, enabling direct application of efficient structured-grid multigrid algorithms. The method preserves second-order accuracy, maintains the optimal complexity of geometric multigrid, and is straightforward to integrate into existing structured-grid solvers.\\
The proposed approach is validated with benchmark problems with analytically known solutions as well as on physically relevant non-uniform geometries. Comparative tests for time efficiency against standard method such as successive over relaxation (SOR) based iterative method and our multigird poisson solver for non uniform grid based on finite difference solvers demonstrate that the transformed-domain finite difference method achieves comparable accuracy along with significantly faster implementation. \\
The structure of the paper is organized as follows: In section 2 we will describe the formulation of method , basic idea of transformation optics and its implementation in computational physics. Section 3 details the geometric multigrid implementation algorithm on uniform computational grid. Section 4 presents validation of numerical solutions obtained using our method with well known analytical solutions and the computational time efficiency with SOR method.

\section{Problem formulation}
\vspace{5mm}
Consider a differentiable coordinate transformation $\mathbf{\mathrm{\textbf{x}\rightarrow \textbf{x}^{\prime}}}$ and  let $\mathcal{J}$ be defined as $\mathrm{3 \times 3}$ Jacobian matrix as
\begin{equation}
\label{eq1}
\mathrm{\mathcal{J}_{ij}=\dfrac{\partial x^{\prime}_{i}}{\partial x_{j}}}.
\end{equation}
The Maxwell's equations can be written in such way  that the  differential form of the equations is invariant under coordinate transforms ~\cite{johnson2007coordinate}. The following vector transformations can be applied to Maxwell's equation so that the form remain invariant under coordinate transformation as 
 \begin{equation}
  \label{eq2}
   \mathrm{\textbf{E}^{\prime} = (\mathcal{J}^{T})^{-1} \mathbf{E}} \ \ \ \ \Rightarrow \ \ \ \ \mathrm{\textbf{E}=\mathcal{J}^{T} \ \textbf{E}^{\prime}}
   \end{equation}
 \begin{equation}
  \label{eq3}
   \mathrm{\textbf{B}^{\prime} = (\mathcal{J}^{T})^{-1} \mathbf{B}} \ \ \ \ \Rightarrow \ \ \ \ \mathrm{\textbf{B}=\mathcal{J}^{T} \ \textbf{B}^{\prime}}
   \end{equation}
  \begin{equation}
  \label{eq4}
   \mathrm{\bm{\nabla}^{\prime} = \bm{\mathcal{J}}^{-1} \bm{\nabla}} \ \ \ \ \Rightarrow \ \ \ \ \mathrm{\bm{\nabla}=\bm{\mathcal{J}} \  \bm{\nabla}^{\prime}}
 \end{equation}
\begin{equation}
  \label{eq5}
       \mathrm{{\rho}^{\prime} = \frac{\rho}{det(\bm{\mathcal{J}})}} \ \ \ \ \Rightarrow \ \ \ \ \mathrm{\rho=det(\bm{\mathcal{J}}) \  \rho^{\prime}}
 \end{equation}
 \begin{equation}
	\label{eq6}
	{\phi}^{\prime} =\phi
\end{equation}
Here, $\mathrm{B^{\prime},E^{\prime},\nabla^{\prime},\rho^{\prime},\phi^{\prime}}$ are  magnetic field, electric field, differential operator, charge density, potential in transformed non uniform domain, while corresponding unprimed quantities are defined over the uniform domain. We will denote the non-uniform (primed) domain as the physical domain and the uniform (unprimed) domain as the computational domain.  \\ 
In this paper, we will only discuss the application regarding electrostatic case involving the Poisson's equation/Gauss's law. Let's consider the equation on the uniform computational grid in the form  
\begin{equation}
    \label{eq7}
    \bm{\nabla}\cdot (\bm{\varepsilon}\textbf{E})=\rho, 
\end{equation}
where $\mathrm{\rho}$ is the free charge density, $\bm{\varepsilon}(\textbf{x})$ is  $3\times3$ relative permittivity tensor and $\mathrm{\textbf{E}}$ is an electric field vector in a uniform computational grid. If we apply vector transformations given in Eqs.~(\ref{eq2})--(\ref{eq6}) in eq (\ref{eq7}) we can simply transform this equation on non-uniform grid as 
\begin{equation}
    \label{eq8}    \nabla^{\prime}\cdot\left(\bm{\varepsilon^{\prime}}\bm{\textbf{E}^{\prime}}\right)=\rho^{\prime},
\end{equation}
provided that the permittivities in both systems transform as 
\begin{equation}
\label{eq9}
 \bm{\varepsilon^{\prime}} = \dfrac{\bm{\mathcal{J}} \bm{\varepsilon} \bm{\mathcal{J}}^{T}}{det(\bm{\mathcal{J}})}
 \end{equation}
 In the physical domain, involving complex geometries, we choose $\mathrm{\mathrm{\varepsilon^{\prime}=I}}$, the identity matrix. Therefore, the permittivity on the computational domain is calculated as 
 \begin{equation}
\label{eq10}
\bm{\varepsilon}=\text{det}(\bm{\mathcal{J}})\bm{\mathcal{J}}^{-1}\bm{\varepsilon}^{\prime}(\bm{\mathcal{J}}^{T})^{-1}. 
\end{equation}
The corresponding charge density on the computation domain is calculated using Eq.(~\ref{eq5}) as 
 \begin{equation}
 \label{eq11}
\rho=\rho^{\prime} det(\bm{\mathcal{J}}).
\end{equation} 
Once the potential and fields are computed on a uniform grid, these quantities can always be transformed back on the physical grid using Eq.(~\ref{eq6}) and (~\ref{eq2}), respectively. 

\subsection{Formulation in two dimensional space}
\vspace{5mm}
For the purpose of demonstration, in this paper we consider only two-dimensional physical grid. The extension to the three-dimensional case is trivial. 
The $\mathrm{2\times2}$ Jacobian matrix, incorporating the coordinate transformation, is written as 
$$\bm{\mathcal{J}_{ij}}=\dfrac{\partial x^{\prime}_{i}}{\partial x_{j}}=\begin{pmatrix}   
	\mathcal{J}_{xx}&\mathcal{J}_{xy}\\[5pt]
	\mathcal{J}_{yx}&\mathcal{J}_{yy}.
\end{pmatrix} \ \ \ $$ 
\\ Using Eq. (\ref{eq10}), we can write the components of the permittivity tensor in a uniform computational space as 
\begin{equation}
\label{eq12}
\mathrm{\varepsilon_{xx}=\dfrac{\mathcal{J}_{yy}^{2}+\mathcal{J}_{xy}^{2}}{det(\mathcal{J})}},
\end{equation}
\begin{equation}
\label{eq13}
\mathrm{\varepsilon_{xy}=\dfrac{-\mathcal{J}_{yy}\mathcal{J}_{yx}-\mathcal{J}_{xy}\mathcal{J}_{xx}}{det(\mathcal{J})}},
\end{equation}
\begin{equation}
\label{eq14}
\mathrm{\varepsilon_{yx}=\dfrac{-\mathcal{J}_{yx}\mathcal{J}_{yy}-\mathcal{J}_{xx}\mathcal{J}_{xy}}{det(\mathcal{J})}},
\end{equation}
\begin{equation}
\label{eq15}
\mathrm{\varepsilon_{yy}=\dfrac{\mathcal{J}_{yx}^{2}+\mathcal{J}_{xx}^{2}}{det(\mathcal{J})}}.
\end{equation} 
Clearly, the modified permittivity on a uniform computational grid is $\mathrm{\varepsilon\neq I}$. \\
We can now write the generalized Poisson equation in the form as\\ $\bm{\nabla}\cdot(\bm{\varepsilon}\textbf{E})=\rho$ \ \ , \ \ \ where $\bm{\varepsilon}=\begin{pmatrix} \varepsilon_{xx}&\varepsilon_{xy}\\ \varepsilon_{yx}&\varepsilon_{yy}\end{pmatrix} \ $ is  $\ 2\times2\ $ modified permittivity tensor. \\ 
This above equation can be expanded as 
\begin{equation}
\label{eq16}
\dfrac{\partial}{\partial x}(\varepsilon_{xx}E_{x})+\dfrac{\partial}{\partial x}(\varepsilon_{xy}E_{y})+\dfrac{\partial}{\partial y}(\varepsilon_{yx}E_{x})+\dfrac{\partial}{\partial y}(\varepsilon_{yy}E_{y})=\rho
\end{equation}
Under electrostatic consideration , putting  $E_{x}=-\dfrac{\partial \phi}{\partial x} \ \ $ and 
		$E_{y}=-\dfrac{\partial \phi}{\partial y} \ \ $ in above equation  , we get 
\begin{equation}
\label{eq17}
\begin{split}
&\varepsilon_{xx}\dfrac{\partial^{2}\phi}{\partial x^{2}}
+ \dfrac{\partial \phi}{\partial x}\dfrac{\partial \varepsilon_{xx}}{\partial x}
+ \varepsilon_{xy} \dfrac{\partial^{2}\phi}{\partial x \partial y} + \dfrac{\partial \phi}{\partial y}\dfrac{\partial \varepsilon_{xy}}{\partial x}
+ \varepsilon_{yx} \dfrac{\partial^{2}\phi}{\partial y \partial x} \\
&+ \dfrac{\partial \phi}{\partial x}\dfrac{\partial \varepsilon_{yx}}{\partial y}
+ \varepsilon_{yy} \dfrac{\partial^{2}\phi}{\partial y^{2}}
+ \dfrac{\partial \phi}{\partial y}\dfrac{\partial \varepsilon_{yy}}{\partial y}
= -\rho
\end{split}
\end{equation}
Using central finite difference second order method we can discretize and rearrange the equation as  
\begin{align}
\label{eq18}
\begin{split}
\phi_{n+1}[i,j] &= 
  k_{1} \phi_{n}[i+1,j] 
 + k_{3} \phi_{n+1}[i-1,j] 
 + k_{4} \phi_{n}[i,j+1] \\
 &+ k_{5} \phi_{n+1}[i,j-1] 
 + k_{6} \big( 
   \phi_{n}[i+1,j+1] - \\
   &\phi_{n+1}[i-1,j+1] 
   - \phi_{n}[i+1,j-1] + \phi_{n+1}[i-1,j-1] 
   \big) \\
 & - \rho[i,j]
 / k_{2}
   \end{split}
\end{align}
Where i,j are the grid indices and 
\begin{align*}
k_{2} &= 2\left(
    \frac{\epsilon_{xx}[i,j]}{\Delta x^{2}}
  + \frac{\epsilon_{yy}[i,j]}{\Delta y^{2}}
\right), \\[6pt]
k_{1} &= \frac{\epsilon_{xx}[i,j]}{\Delta x^{2}}
  + \frac{
      \epsilon_{xx}[i+1,j]-\epsilon_{xx}[i-1,j]
    }{4\Delta x^{2}} \notag \\
&\quad + \frac{
      \epsilon_{yx}[i,j+1]-\epsilon_{yx}[i,j-1]
    }{4\Delta x \Delta y}, \\[6pt]
k_{3} &= \frac{\epsilon_{xx}[i,j]}{\Delta x^{2}}
  - \frac{
      \epsilon_{xx}[i+1,j]-\epsilon_{xx}[i-1,j]
    }{4\Delta x^{2}} \notag \\
&\quad - \frac{
      \epsilon_{yx}[i,j+1]-\epsilon_{yx}[i,j-1]
    }{4\Delta x \Delta y}, \\[6pt]
k_{4} &= \frac{\epsilon_{yy}[i,j]}{\Delta y^{2}}
  + \frac{
      \epsilon_{yy}[i,j+1]-\epsilon_{yy}[i,j-1]
    }{4\Delta y^{2}} \notag \\
&\quad + \frac{
      \epsilon_{xy}[i+1,j]-\epsilon_{xy}[i-1,j]
    }{4\Delta x \Delta y}, \\[6pt]
k_{5} &= \frac{\epsilon_{yy}[i,j]}{\Delta y^{2}}
  - \frac{
      \epsilon_{yy}[i,j+1]-\epsilon_{yy}[i,j-1]
    }{4\Delta y^{2}} \notag \\
&\quad - \frac{
      \epsilon_{xy}[i+1,j]-\epsilon_{xy}[i-1,j]
    }{4\Delta x \Delta y}, \\[6pt]
k_{6} &= \frac{\epsilon_{xy}[i,j]}{2\Delta x \Delta y}.
\end{align*}
Then we can write residue as $$b[i,j]=\Phi_{n+1}[i,j]-\Phi_{n}[i,j]$$
and the final solution can be written as $$\Phi_{n+1}[i,j]=\Phi_{n}[i,j]+\omega\times b[i,j]$$
Here $\omega$ is a relaxation factor and $1<\omega<2$ . If for $\omega=1$, the method reduces to Gauss-Seidel method, which is used for smoothening of the solution in the multigrid algorithm. The implementation of multi-grid algorithm will be discussed in detail in the next section. \\
In summary, the key idea of this algorithm is to apply coordinate transformation $\mathbf{\mathrm{\textbf{x}\rightarrow \textbf{x}^{\prime}}}$ to get the required non-uniform grid ($\textbf{x}^{\prime}$) from the uniform computational grid ($\textbf{x}$). The corresponding Jacobian matrix, calculated from Eq.(~\ref{eq1}), is used to computed the components of permittivity tensor ($\epsilon$), and charge densities ($\rho$) using Eq.(\ref{eq10}) and (\ref{eq11}), respectively. With modified $\epsilon$ and $\rho$ on the uniform grid, we solve generalized Poisson equation with finite difference method. Finally, the solution obtained on the uniform grid is transformed back on the physical grid using the inverse transformation.
\section{Implementation with Geometric Multigrid Method}
\vspace{5mm}
The Geometric Multigrid (GMG) method is a well-established iterative technique for solving large, sparse systems of equations arising from the discretization of partial differential equations. GMG operates on a hierarchy of grids generated through systematic coarsening of the original computational mesh, transferring information between fine and coarse levels via restriction and prolongation operators. This approach is highly efficient for structured meshes and problems with simple, well-defined geometries, where coarser grids can be constructed naturally by halving the mesh spacing. However, GMG’s reliance on a regular geometric hierarchy limits its applicability on irregular domains, where generating consistent coarser levels may require complex remeshing and can lead to loss of accuracy near curved or highly nonuniform boundaries.\\
In contrast, the Algebraic Multigrid (AMG) method constructs its grid hierarchy directly from the system matrix, avoiding the need for explicit geometric information. This makes AMG naturally suited for unstructured meshes, complex geometries, and strongly heterogeneous coefficients. Nevertheless, AMG can experience performance degradation on irregular domains if the interpolation operators fail to accurately capture the smooth error components—particularly in the presence of strong anisotropy, abrupt coefficient jumps, or distorted mesh cells. In such situations, advanced AMG variants such as smoothed aggregation or adaptive AMG are employed to maintain convergence efficiency. While AMG generally incurs a higher per-iteration computational cost than GMG, its robustness and flexibility often make it the preferred choice for problems involving complex geometrical configurations.\\
In the present work, the algorithm described in the previous section is ideally suitable for the GMG methods, as all the computations are performed on a uniform grid. Therefore, once the components of permittivity tensors on the structured grid are pre-calculated, the additional overhead of interpolation, typically encountered in AMG, is eliminated. This leads to straight-forward implementation of GMG for problems which conventionally require non-uniform grid. Additionally, the structured uniform mesh allows easy parallelization of the problem. The implementation of GMG is briefly described below.\\
Consider a discrete system resulting from the discretization of a continuous problem on a regular(uniform) mesh with grid spacing \( h \):
\begin{equation}
    \label{eq19}
    \mathbf{A}_h \Phi_h = \mathbf{b}_h,
\end{equation}
where \( \mathbf{A}_h \) is the discretized operator, \( \Phi_h \) is the unknown solution vector, and \( \mathbf{b}_h \) is the source term.
Let \( \bar{\Phi}_h \) be the current approximate solution obtained via an iterative method. Since \( \bar{\Phi}_h \) is not exact, it produces a residual
\begin{equation}
    \label{eq20}
    \mathbf{r}_h = \mathbf{b}_h - \mathbf{A}_h \bar{\Phi}_h,
\end{equation}
which vanishes only when the exact solution is found.
The goal of the multigrid method is to compute a correction \( \Psi_h \) such that
\begin{equation}
    \label{eq21}
    \Phi_h = \bar{\Phi}_h + \Psi_h.
\end{equation}
Substituting this into \eqref{eq19} and using linearity of \( \mathbf{A}_h \) yields
\begin{equation}
    \label{eq22}
    \mathbf{A}_h \Psi_h = \mathbf{r}_h.
\end{equation}
High-frequency errors in \( \bar{\Phi}_h \) are reduced by smoothing iterations on the fine grid. The remaining smooth error \( \Psi_h \) can be more efficiently approximated on a coarser grid with mesh size \( H > h \) by solving
\begin{equation}
    \label{eq23}
    \mathbf{A}_H \Psi_H = \mathbf{r}_H,
\end{equation}
where the coarse-grid residual is obtained by restricting the fine-grid residual:
\[
\mathbf{r}_H = \mathbf{I}_h^H \mathbf{r}_h,
\]
with \( \mathbf{I}_h^H \) denoting the restriction operator.
The coarse-grid correction \( \Psi_H \) is then interpolated (prolongated) back to the fine grid to update the solution:
\begin{equation}
    \label{eq24}
    \bar{\Phi}_h^{\text{new}} = \bar{\Phi}_h + \mathbf{I}_H^h \Psi_H,
\end{equation}
where \( \mathbf{I}_H^h \) is the prolongation operator.
Additional smoothing steps on the fine grid reduce high-frequency errors introduced by interpolation. This sequence of operations forms one multigrid cycle, repeated until convergence. 
\section{Results and Benchmarking}
\vspace{5mm}
For benchmarking, we validate our results against well-known analytical solutions and compare computational times between the GMG method and the Successive Over-Relaxation (SOR) iterative method. The results are classified into two parts: (i) grids generated via analytical transformation, and (ii) grids generated via numerical transformation. For direct SOR simulations, we use $\omega = 1.875$ and a tolerance of $10^{-10}$. \\
For simplicity, we choose an equal number of points in both $x$ and $y$ directions, i.e., $\mathrm{N_x = N_y = N}$. The error analysis is performed by computing the following error norms over the 2-D domain:
\begin{equation}
\label{eq25}
    \mathrm{L_{\infty} = \underset{(x_i,y_j) \in D}{\mathrm{max}} \ \left| \phi_{\mathrm{exact}}(x_i,y_j) - \phi_{\mathrm{numerical}}(x_i,y_j) \right| ,}
\end{equation}
\begin{equation}
\label{eq26}
    \mathrm{L_{2} = \sqrt{\frac{1}{\mathrm{N}^2} \sum_{(x_i,y_j) \in D} \left| \phi_{\mathrm{exact}}(\mathrm{x_i,y_j}) - \phi_{\mathrm{numerical}}(\mathrm{x_i,y_j}) \right|^2} ,}
\end{equation}
\begin{equation}
\label{eq27}
    \mathrm{L_{1} = \frac{1}{N^2} \sum_{(x_i,y_j) \in D} \left| \phi_{\mathrm{exact}}(x_i,y_j) - \phi_{\mathrm{numerical}}(x_i,y_j) \right| .}
\end{equation}
\subsection{Result through analytical grid generation}
\begin{figure}[H]
	\centering
	\includegraphics[width=1.0\textwidth]{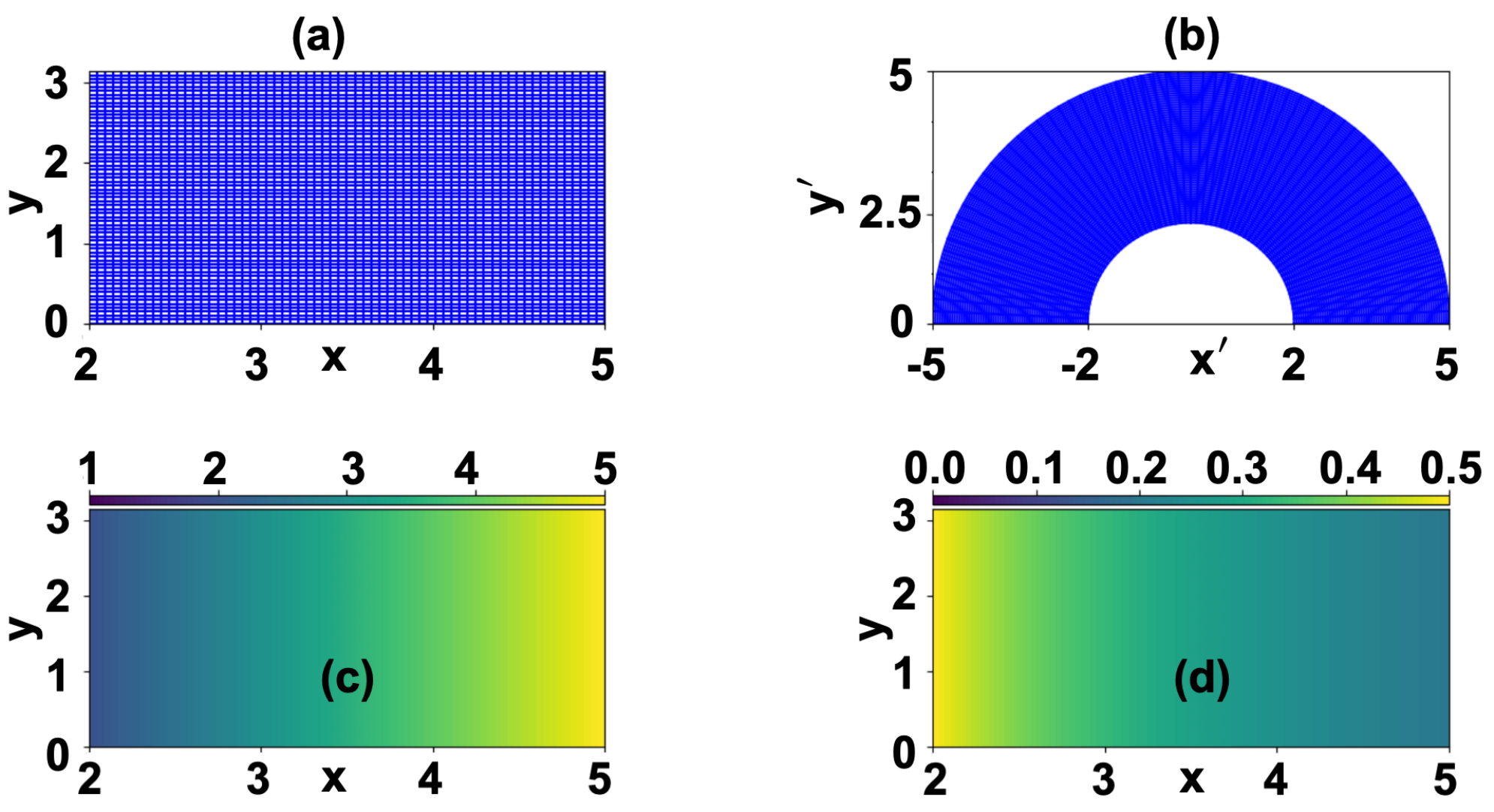}
	\caption{(a) Computational grid, (b) Physical grid , (c) $\varepsilon_{xx}$ , (d) $\varepsilon_{yy}$ . }
	\label{cylindrical_grid_eps}
    \end{figure}

\begin{figure}[H]
    \centering
	\includegraphics[width=1.0\textwidth]{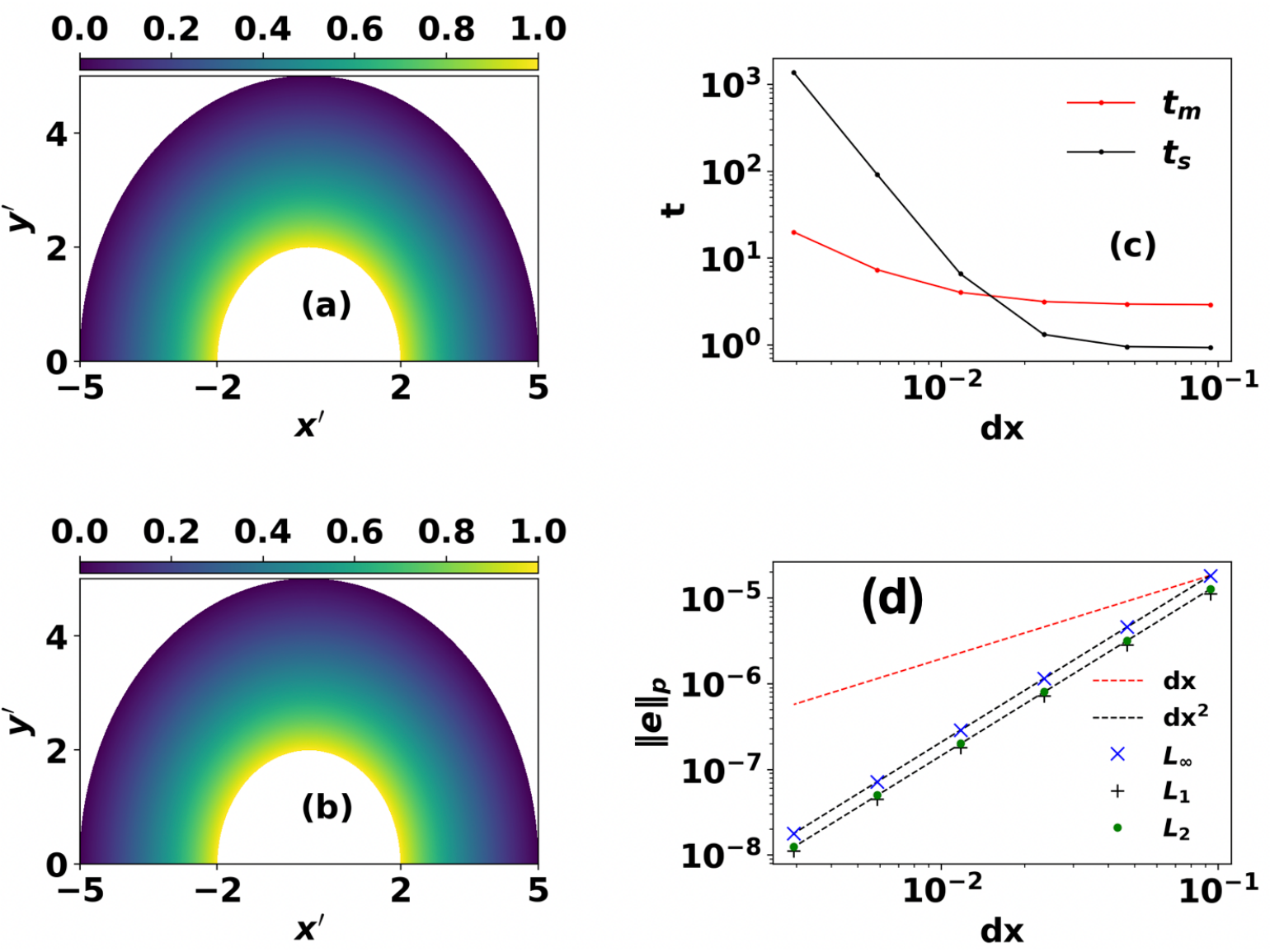}
	\caption{(a) Analytical solution, (b) Numerical solution, (c) Speed-up comparison of Multigrid and SOR method, (d) error analysis. The second order accuracy of the scheme can be seen from the subplot (d).}
	\label{cylindrical_solution}
\end{figure}

\vspace{3mm}
For certain problems in physics and engineering, the analytical mapping between the computational grid to physical grid is obvious. In this situation, the exact expressions of the permittivity tensor are known from this mapping. For the purpose of demonstration, we consider a semi-annulus domain shown in Fig.~\ref{cylindrical_grid_eps}. The computational grid and the physical grid are shown in Fig.~\ref{cylindrical_grid_eps} (a) and (b), respectively. On the inner wall of this cylinder ($r = a$), a finite potential $\Phi_0$ is applied. The outer wall ($r = R$) is grounded with the region between the two walls is taken as charge-free ($\mathrm{\rho^{\prime} = 0}$). \\    
These boundary conditions are expressed as follows: 
\begin{equation}
\label{eq28}
	\begin{split}
		\dfrac{\mathrm{\partial \Phi}}{\mathrm{\partial n}}&=0 \ \ \ \ \   \text{at} \ \ \theta =0 \\[8pt]
		\dfrac{\mathrm{\partial \Phi}}{\mathrm{\partial n}}&=0 \ \ \ \ \   \text{at} \ \ \theta =\pi \\[8pt]
		\Phi &= 0  \ \ \ \ \ \text{at} \ \ \mathrm{r = R} \\[8pt]  
		\Phi &= \Phi_{0}  \ \ \ \text{at} \ \ \mathrm{r = a} \\[8pt]
	\end{split} 
\end{equation}
The exact solution, $\mathrm{\Phi(r,\theta)}$, for this problem can be written as:
\begin{equation}
    \label{eq28-a}
    \mathrm{\Phi(r,\theta)} =  \dfrac{\Phi_{0} \ln{\mathrm{\left(\dfrac{R}{r}\right)}}}{\ln{\mathrm{\left(\dfrac{R}{a}\right)}}}
\end{equation}
The mapping to generate the semi-annulus domain (Fig.\ref{cylindrical_grid_eps}(b)) from the uniform computational domain (Fig.\ref{cylindrical_grid_eps}(a)) is described by a coordinate transformation ~\cite{isshiki2017generation}:
\begin{equation}
\label{eq29}
	\begin{split}
		\mathrm{x^{\prime}} = \mathrm{x \cos{y}} \\[1pt]
		\mathrm{y^{\prime}} = \mathrm{x \sin{y}}
	\end{split}	
\end{equation}  
We set the inner and outer radii to $\mathrm{a = 2}$ and $\mathrm{R=5}$, respectively, and use a grid resolution of $\mathrm{n_x = n_y = 201}$, with corresponding grid spacings of $\mathrm{\Delta x = \left(\dfrac{R - a}{n_x - 1}\right)}$ and $\mathrm{\Delta y = \left(\dfrac{\pi}{n_y - 1}\right)}$. The potential on the inner boundary is set to $\mathrm{\Phi_0 = 1}$. The diagonal components of the permittivity tensor, viz. $\mathrm{\varepsilon_{xx}}$ and $\mathrm{\varepsilon_{yy}}$, obtained from Eq.~\ref{eq12}-\ref{eq15}, are shown in Fig.~\ref{cylindrical_grid_eps}(c) and (d), respectively. Note that the off-diagonal components vanish for this transformation i.e. $\mathrm{\varepsilon_{xy}=\varepsilon_{yx}=\mathbf{0}}$. \\  
The analytical solution is compared with the numerical solution in Figure~\ref{cylindrical_solution} (a) and (b), respectively. The error analysis, presented in Fig.~\ref{cylindrical_solution}(d), demonstrates second-order accuracy in both local ($\mathrm{L^{\infty}}$) and global ($\mathrm{L^1}$, $\mathrm{L^2}$) norms. The reference dotted lines are drawn for the dependence of errors with the first ($\propto \mathrm{\Delta x}$) and second ($ \propto \mathrm{\Delta x^2}$) order accuracies. Here, we see that all three norms fall on the second order accuracy line.  In addition, Fig.~\ref{cylindrical_solution}(c) illustrates that the multigrid solver achieves a speed-up of approximately two orders of magnitude compared to the SOR method. Thus, the algorithm achieves significant speed-up with minimal changes in the standard textbook Poisson solvers.  Note that all runtime reported in the paper are obtained by running the code on laptop machine with Intel Core i7-1165G7 with 16 GB RAM.
\subsection{Result through numerical grid generation}
\vspace{3mm}
 When obvious analytical mapping, for transforming the physical domain from the computational domain, is absent; the generation of physical grid-space through numerical methods is essential. In this section, we demonstrate the extension of algorithm to examples where numerical grid generation is required. Numerical grid generation is widely used in many fields such as computational fluid dynamics~\cite{sanmiguel2005cartesian,eiseman1985grid,thompson1985numerical} and conformal mapping in electromagnetism ~\cite{landy2009guiding}. This approach allows us to deal with complex geometries to be handled by mapping an arbitrary grid to Cartesian grid. The most common approach is based on solving Laplace's equation ~\cite{thompson1985numerical,akcelik2001nearly} for the transformed coordinates because it generates coordinates that vary linearly. Given the boundary conditions for an arbitrary coordinate transformation from $\mathrm{(x,y,z)}$ to $\mathrm{(x^{\prime},y^{\prime},z^{\prime})}$, the Laplace's equation can be used to solve for transformed coordinates using the following equations.
 \begin{equation}
 	\label{eq30} 	\mathrm{\nabla^{2}x^{\prime}=\left(\dfrac{\partial^{2}}{\partial x^{2}}+\dfrac{\partial^{2}}{\partial y^{2}}+\dfrac{\partial^{2}}{\partial z^{2}}\right)x^{\prime} = 0}
 \end{equation}
 \begin{equation}
 	\label{eq31}	\mathrm{\nabla^{2}y^{\prime}=\left(\dfrac{\partial^{2}}{\partial x^{2}}+\dfrac{\partial^{2}}{\partial y^{2}}+\dfrac{\partial^{2}}{\partial z^{2}}\right)y^{\prime} = 0}
 \end{equation}
 \begin{equation}
 	\label{eq32}	\mathrm{\nabla^{2}z^{\prime}=\left(\dfrac{\partial^{2}}{\partial x^{2}}+\dfrac{\partial^{2}}{\partial y^{2}}+\dfrac{\partial^{2}}{\partial z^{2}}\right)z^{\prime} = 0}
 \end{equation}
 \\ In some cases backward transformation Laplace's equation  can be used as 
 \begin{equation}
 	\label{eq33} 	\mathrm{\nabla^{2}x=\left(\dfrac{\partial^{2}}{\partial {x^{\prime}}^{2}}+\dfrac{\partial^{2}}{\partial {y^{\prime}}^{2}}+\dfrac{\partial^{2}}{\partial {z^{\prime}}^{2}}\right)x = 0}
 \end{equation}
 \begin{equation}
 	\label{eq34}	\mathrm{\nabla^{2}y=\left(\dfrac{\partial^{2}}{\partial {x^{\prime}}^{2}}+\dfrac{\partial^{2}}{\partial {y^{\prime}}^{2}}+\dfrac{\partial^{2}}{\partial {z^{\prime}}^{2}}\right)y = 0}
 \end{equation}
 \begin{equation}
 	\label{eq35} 	\mathrm{\nabla^{2}z=\left(\dfrac{\partial^{2}}{\partial {x^{\prime}}^{2}}+\dfrac{\partial^{2}}{\partial {y^{\prime}}^{2}}+\dfrac{\partial^{2}}{\partial {z^{\prime}}^{2}}\right)z = 0}
 \end{equation}
For the 2-D case, we can solve for transformed system as 
 \begin{equation}
 	\label{eq36} 	\mathrm{\nabla^{2}x^{\prime}=\left(\dfrac{\partial^{2}}{\partial x^{2}}+\dfrac{\partial^{2}}{\partial y^{2}}\right)x^{\prime} = 0},
 \end{equation}
 \begin{equation}
 	\label{eq37} 	\mathrm{\nabla^{2}y^{\prime}=\left(\dfrac{\partial^{2}}{\partial x^{2}}+\dfrac{\partial^{2}}{\partial y^{2}}\right)y^{\prime} = 0},
 \end{equation}
 subject to appropriate boundary condition of the form $\mathrm{x^{\prime}=f(x,y) \ and \ \ y^{\prime}=f(x,y)}$. We can now compute the coordinates of the numerically generated grid by using  Eqs.~(\ref{eq36}), (\ref{eq37}) by applying the Geometric Multigrid algorithm  as described in section 3.
The elements of Jacobian matrix can be pre-calculated numerically using the second order finite difference formula as follows: 
\begin{equation}
	\label{eq39}
\mathrm{\mathcal{J}_{xx}[i,j] = \dfrac{\partial x^{\prime}}{\partial x} \cong \dfrac{x^{\prime}(i+1,j)-x^{\prime}(i-1,j)}{2\Delta x}}, 
\end{equation}
\\ \begin{equation}
	\label{eq40}
\mathrm{\mathcal{J}_{xy}[i,j] = \dfrac{\partial x^{\prime}}{\partial y} \cong \dfrac{x^{\prime}(i,j+1)-x^{\prime}(i,j-1)}{2\Delta y}}, 
\end{equation}
\\ \begin{equation}
	\label{eq41}
	\mathrm{\mathcal{J}_{yx}[i,j] = \dfrac{\partial y^{\prime}}{\partial x} \cong \dfrac{y^{\prime}(i+1,j)-y^{\prime}(i-1,j)}{2\Delta x}}, 
\end{equation}
\\ \begin{equation}
	\label{eq42}
	\mathrm{\mathcal{J}_{yy}[i,j] = \dfrac{\partial y^{\prime}}{\partial y} \cong \dfrac{y^{\prime}(i,j+1)-y^{\prime}(i,j-1)}{2\Delta y}}. 
\end{equation}
\\ These components of the Jacobian matrix are substituted in Eqs.~(\ref{eq12})-(\ref{eq15}) for the computation of the permittivity tensor at any grid location $\mathrm{(i,j)}$. The method is demonstrated through specific examples discussed below.

\subsubsection{Circular Domain with a Source Term}
\hfill\break

\begin{figure}[H]
    \centering
    \includegraphics[width=1.0\textwidth]{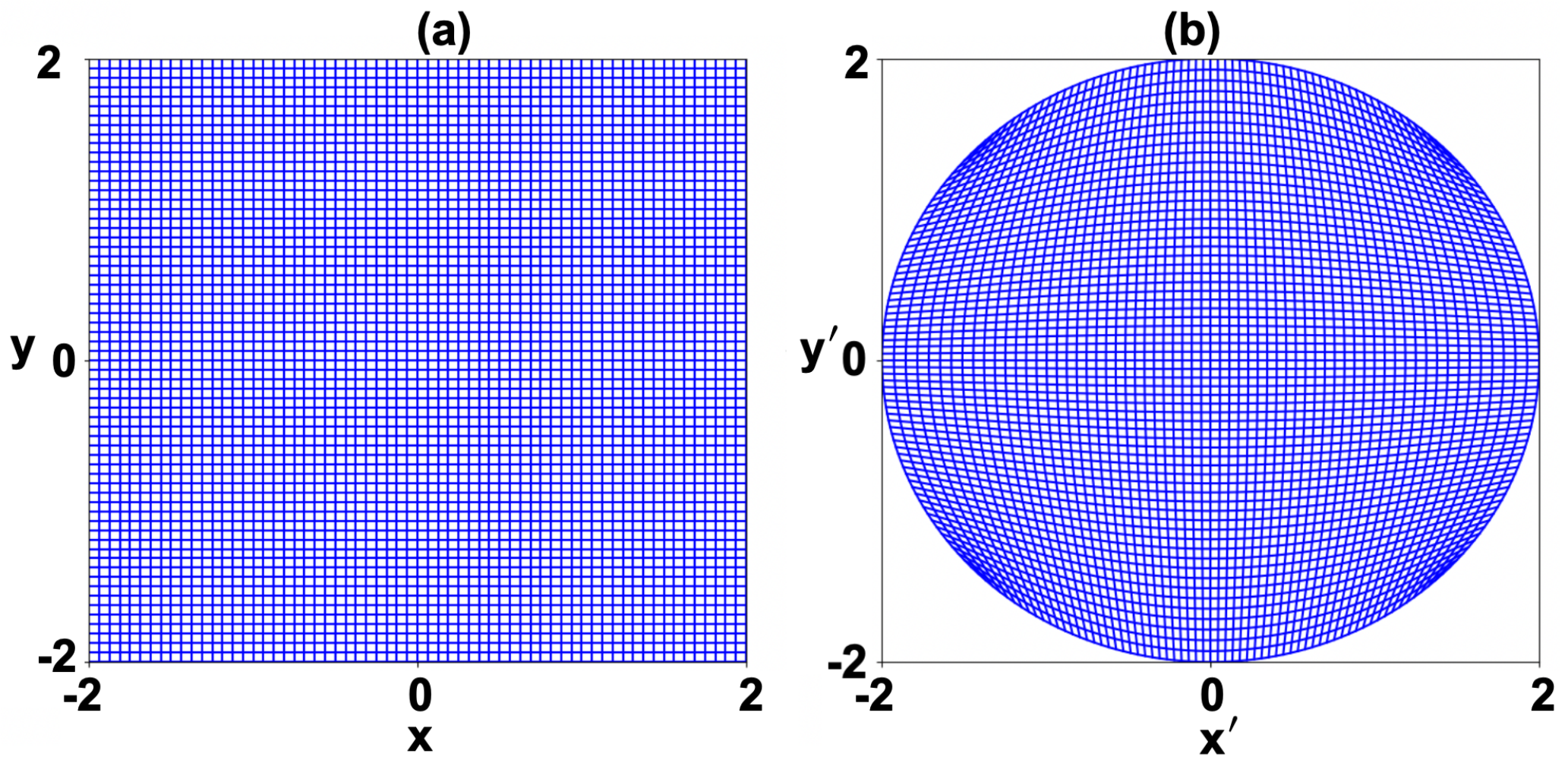}
    \caption{Circular domain mapping: (a) Uniform computational grid and (b) transformed physical grid.}
    \label{num_cir_grid}
\end{figure}

\begin{figure}[H]
    \centering
    \includegraphics[width=1.0\textwidth]{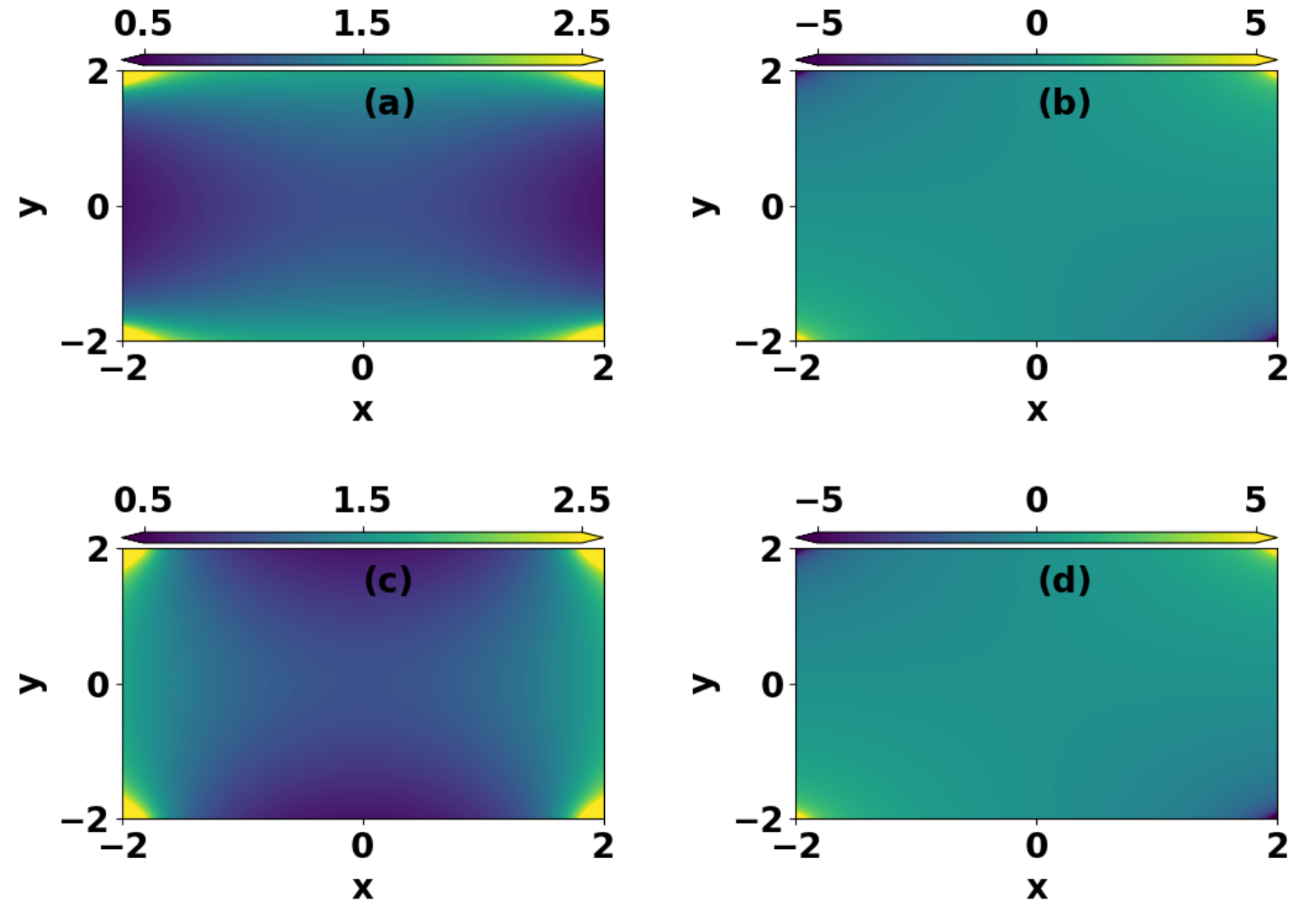}
    \caption{Components of the permittivity tensor: (a) $\varepsilon_{xx}$, (b) $\varepsilon_{xy}$, (c) $\varepsilon_{yy}$, and (d) $\varepsilon_{yx}$.}
    \label{num_cir_eps}
\end{figure}

\begin{figure}[H]
    \centering
\includegraphics[width=1.0\textwidth]{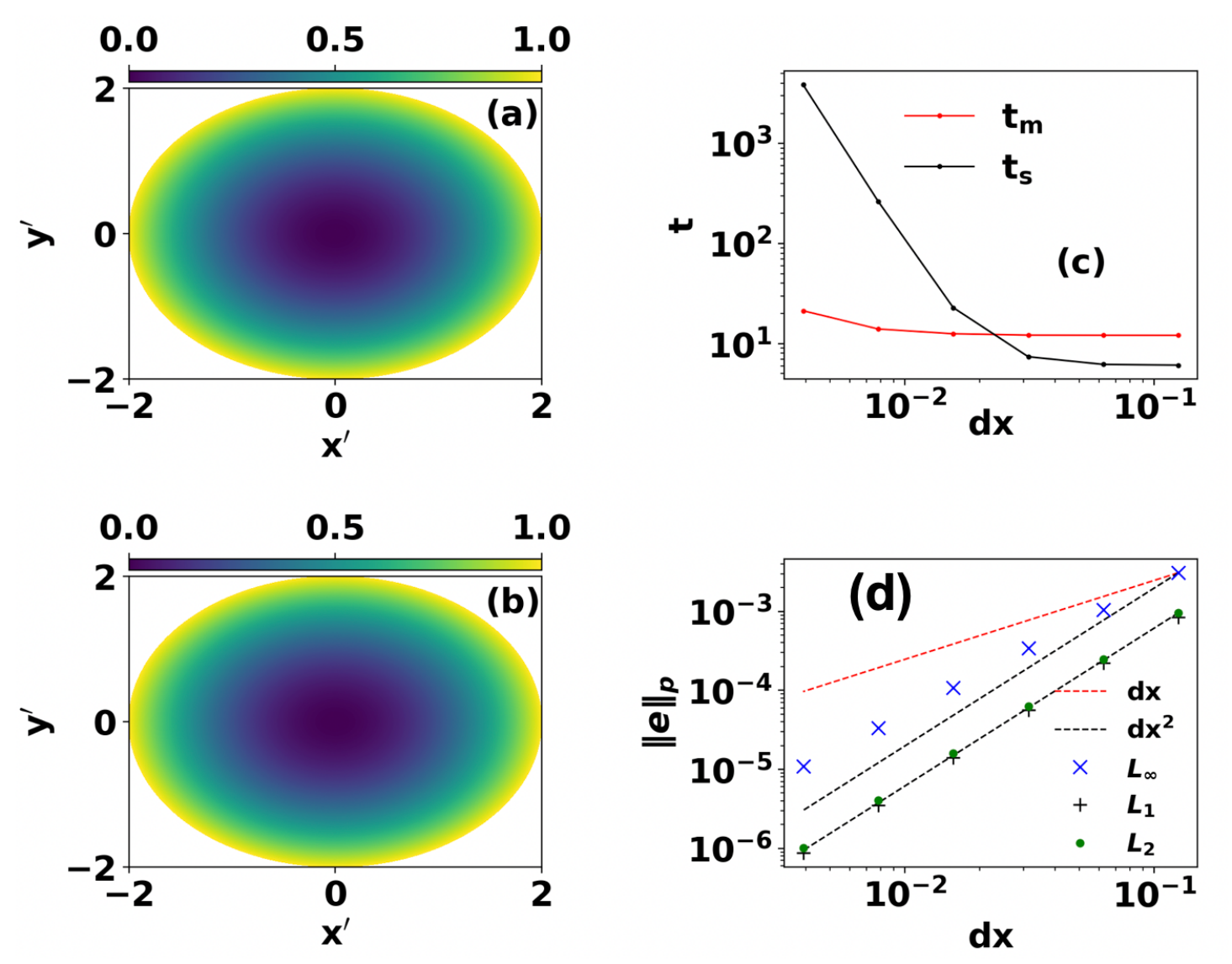}
    \caption{Results for the circular domain: (a) Analytical solution, (b) Numerical solution using GMG, (c) Runtime comparison of SOR and GMG methods, and (d) Grid convergence in $\mathrm{L_1}$, $\mathrm{L_2}$, and $\mathrm{L_{\infty}}$ norms.}
    \label{num_cir_sol}
\end{figure}

In this example, we consider the solution of Poisson equation over a circular domain (see Fig.~\ref{num_cir_grid}(b)) of radius $b$ with uniform source term $\rho^{\prime} = 1$. The value of the potential on the boundary is taken as $b^2/4$. The analytical solution for this problem is expressed as 
\begin{equation}
    \label{eq43}
    \phi_{\mathrm{exact}}(x^{\prime},y^{\prime}) = \frac{{x^{\prime}}^2+{y^{\prime}}^2}{4}
\end{equation}\\
In order to map, rectangular computational domain to circular physical domain, we consider the following transformation:
\begin{equation}
\label{eq44}
    \begin{split}
        &x^{\prime} = \frac{x}{\sqrt{2}}, \quad y^{\prime} = -b\sqrt{1 - \frac{x^{2}}{2b^{2}}} \quad (-b \leq x \leq b \, ; \, y = -b), \\[4pt]
        &x^{\prime} = \frac{x}{\sqrt{2}}, \quad y^{\prime} = \ \ b\sqrt{1 - \frac{x^{2}}{2b^{2}}} \quad (-b \leq x \leq b \, ; \, y = b), \\[4pt]
        &x^{\prime} = -b\sqrt{1 - \frac{y^{2}}{2b^{2}}}, \quad y^{\prime} = \frac{y}{\sqrt{2}} \quad ( x = -b \, ; \, -b \leq y \leq b ), \\[4pt]
        &x^{\prime} = \ \ b\sqrt{1 - \frac{y^{2}}{2b^{2}}}, \quad y^{\prime} = \frac{y}{\sqrt{2}} \quad ( x = b \, ; \, -b \leq y \leq b ) .
    \end{split}
\end{equation}
Note that this transformation maps sides of a square to each quadrant of the circle. Using these conditions as boundary conditions for Eq.~\ref{eq36} and \ref{eq37}, we obtain the physical grid from the cartesian physical grid as shown in Fig.~\ref{num_cir_grid}. Here, we have taken $b = 2$ and the grid resolution is set to $\mathrm{N_x = N_y = 257}$, corresponding to a uniform spacing $\Delta x = \Delta y = 0.015625$. The uniform source term in the physical domain turns into a non-uniform source with the Eq.~\ref{eq5}.
The numerically computed components of the permittivity tensor are shown in Fig.~\ref{num_cir_eps}. 
The comparison between the numerical and analytical solutions is presented in Fig.~\ref{num_cir_sol} (a) and (b), respectively. The grid resolution scan, shown in subplot~(d), confirms second-order accuracy in $\mathrm{L_1}$, $\mathrm{L_2}$, and $\mathrm{L_{\infty}}$ norms. Additionally, Subplot~(c) highlights the computational efficiency of the geometric multigrid method, showing it to be approximately three times faster than the general SOR method.
\subsubsection{Solution of Laplace equation over circular domain }
\hfill\break

\begin{figure}[H]
    \centering
    \includegraphics[width=1.0\textwidth]{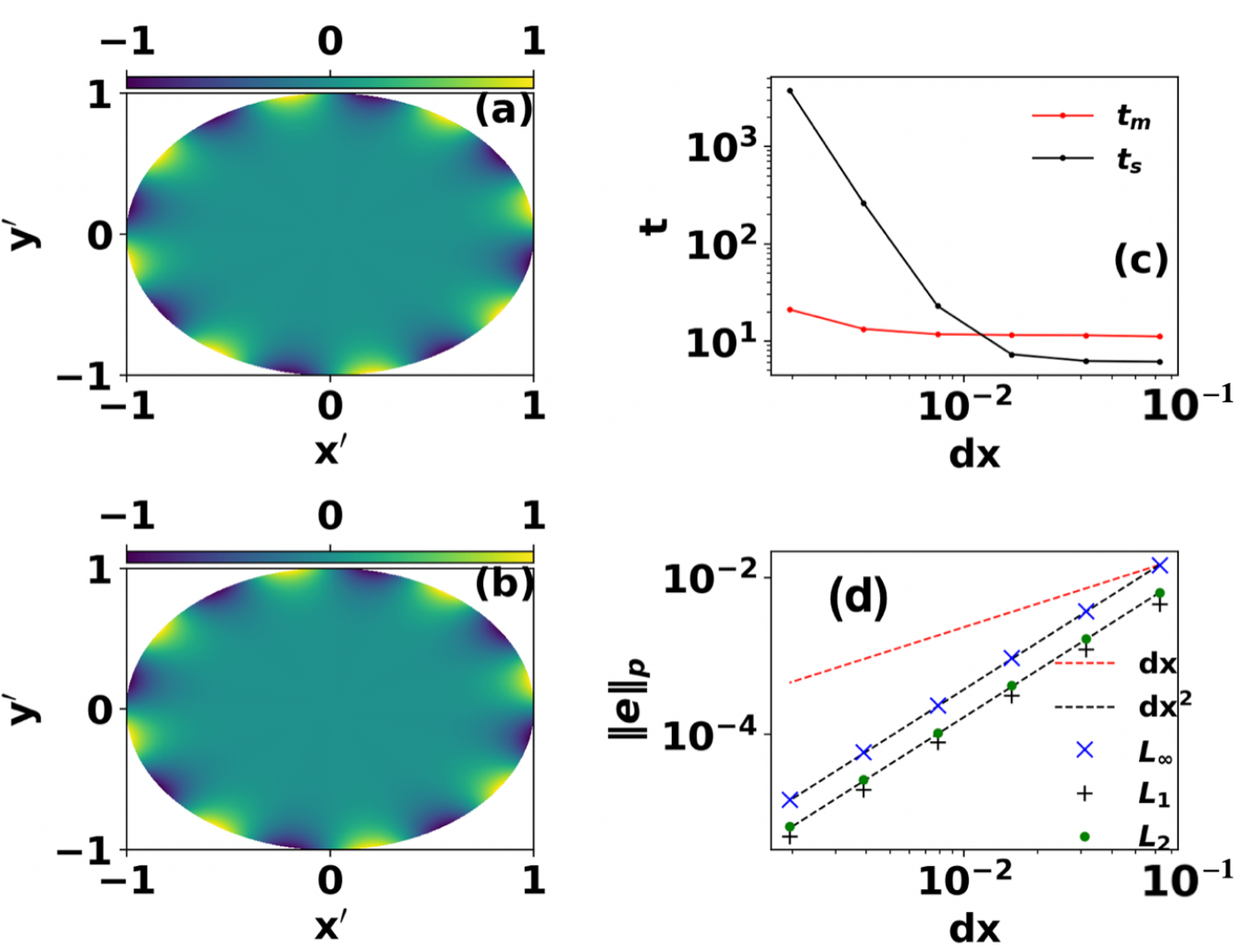}
    \caption{Results for the unit-radius circle: (a) Analytical solution, (b) Numerical solution using GMG, (c) Runtime comparison between SOR and GMG methods, and (d) Grid convergence in $\mathrm{L_1}$, $\mathrm{L_2}$, and $\mathrm{L_{\infty}}$ norms.}
    \label{hm_solution}
\end{figure}

We now consider a solution of Laplace equation ($\rho^{\prime} = 0$) over a unit circle ($b = 1$). We apply Dirichlet boundary condition in the form
\begin{equation}
    \label{eq45-1}
    \phi(b,\theta) = \sin(n\theta).
\end{equation}
The analytical solution inside the domain ($r \leq 1$) can be written in the form
\begin{equation}
    \label{eq45}
    \phi_{\mathrm{exact}}(r,\theta) = r^n \sin(n \theta)
\end{equation}
In this study, we choose $n = 8$. Using identical mapping given in Eq.~\ref{eq44}, we generate physical and computational grids similar to Fig.~\ref{num_cir_grid}. The size of the computational domain is chosen to $\mathrm{L_x = L_y = 2}$ with a resolution of $\mathrm{N_x = N_y = 257}$. Figure~\ref{hm_solution} show the analytical expression for the potential (subplot a) matches well with the numerical (subplot b) results with GMG method. The error dependence on grid resolution (subplot d) clearly demonstrates second order global and local accuracy. Here, it should be noted that identical problem was solved with the boundary element method (BEM) by Rapaka et al ~\cite{rapaka2020efficient}. The solution with BEM method showed first order accuracy in local ($\mathrm{L_{\infty}}$) norm and approximately second order accuracy in the global ($\mathrm{L_1}$, $\mathrm{L_2}$) norms (Fig. 6 of Rapaka et al \cite{rapaka2020efficient}). In contract, our method has demonstrated exact second order accuracy in both local and global norms. Additionally, Subfigure~(c) demonstrates the improved computational efficiency of the GMG method, which is approximately three times faster than the classical SOR method. 

\section{Some applications based on the method}
\vspace{5mm}
Having demonstrated benchmarking of the algorithm with standard problems having well-known analytical solutions, we now show the application of this method to some of the problems in plasma physics and fluid dynamics. 

\subsection{Stretched/non-uniform grid simulations for enhanced computational efficiency }
\begin{figure}[H]
	\centering
	\includegraphics[width=1.0\textwidth]{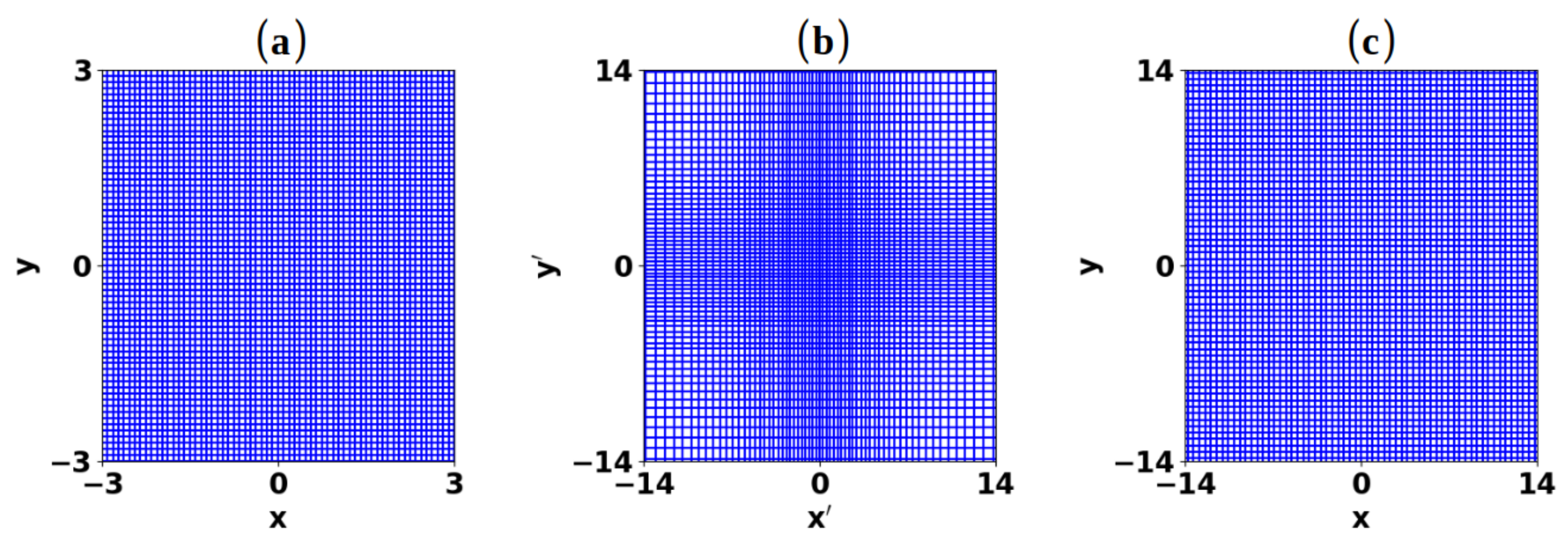}
	\caption{(a) Initial uniform computational grid. (b) Transformed non-uniform physical grid using Eq.~\ref{eq52}. (c) Equivalent uniform grid resolution required to match smallest spacing of the non-uniform grid.}
	\label{stretch_grid}
\end{figure}

\begin{figure}[H]
	\centering
	\includegraphics[width=1.0\textwidth]{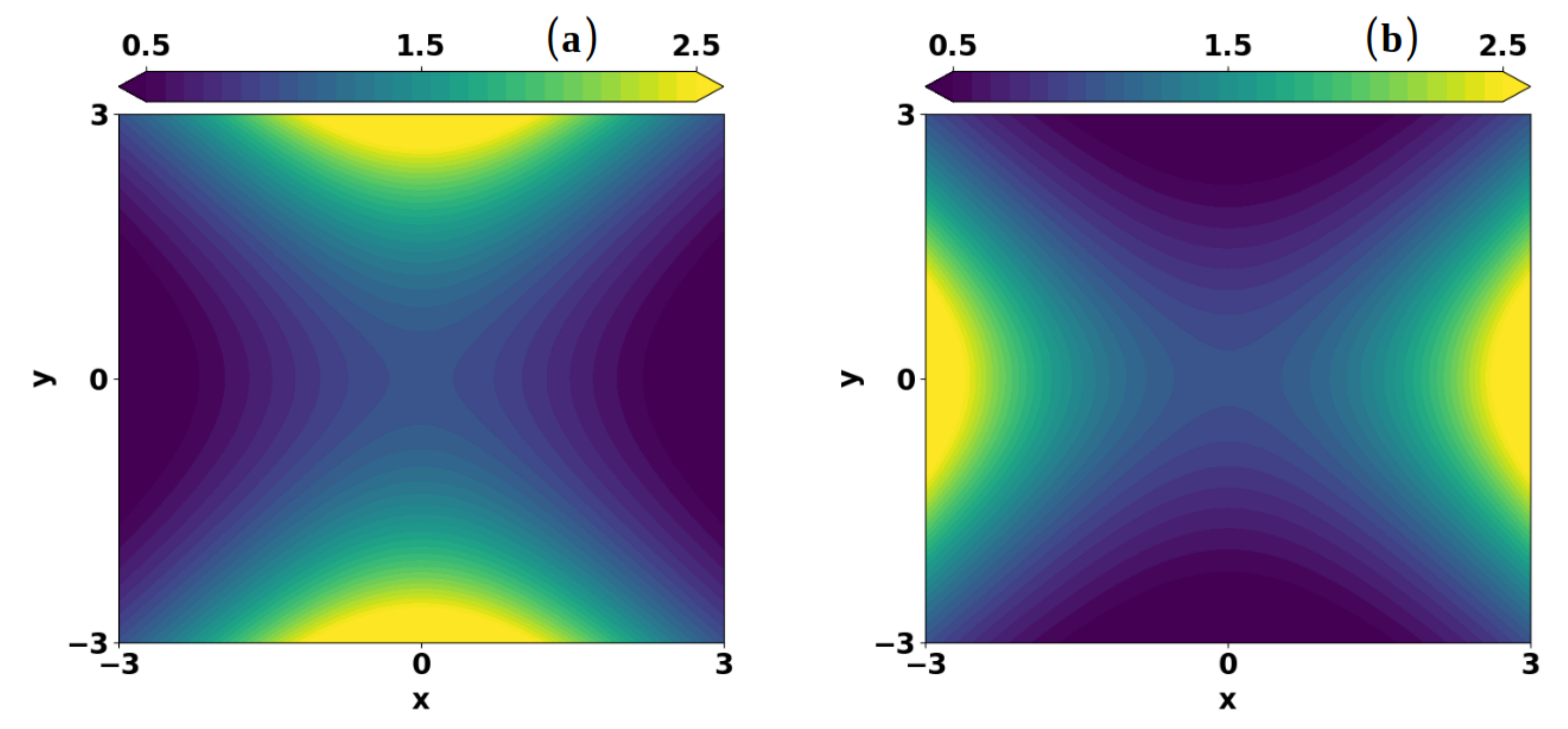}
	\caption{Diagonal components of the anisotropic permittivity tensor obtained from the coordinate transformation.}
	\label{epsilon_stretch}
\end{figure}

\begin{figure}[H]
	\centering
	\includegraphics[width=1.0\textwidth]{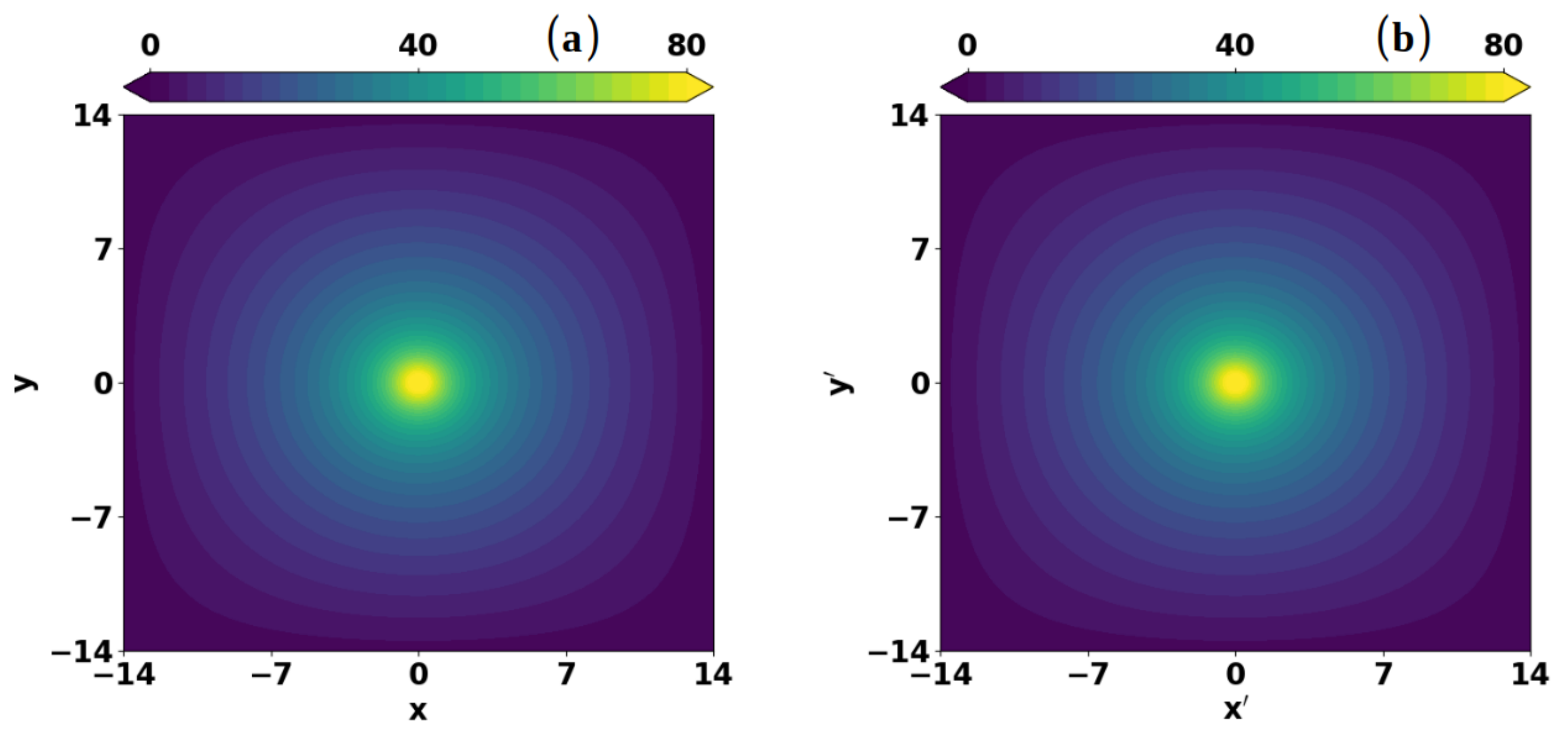}
	\caption{(a) Solution of the Poisson equation on the stretched non-uniform grid. (b) Equivalent solution obtained on a larger uniform grid.}
	\label{potential_stretch}
\end{figure}

\begin{figure}[H]
	\centering
	\includegraphics[width=1.0\textwidth]{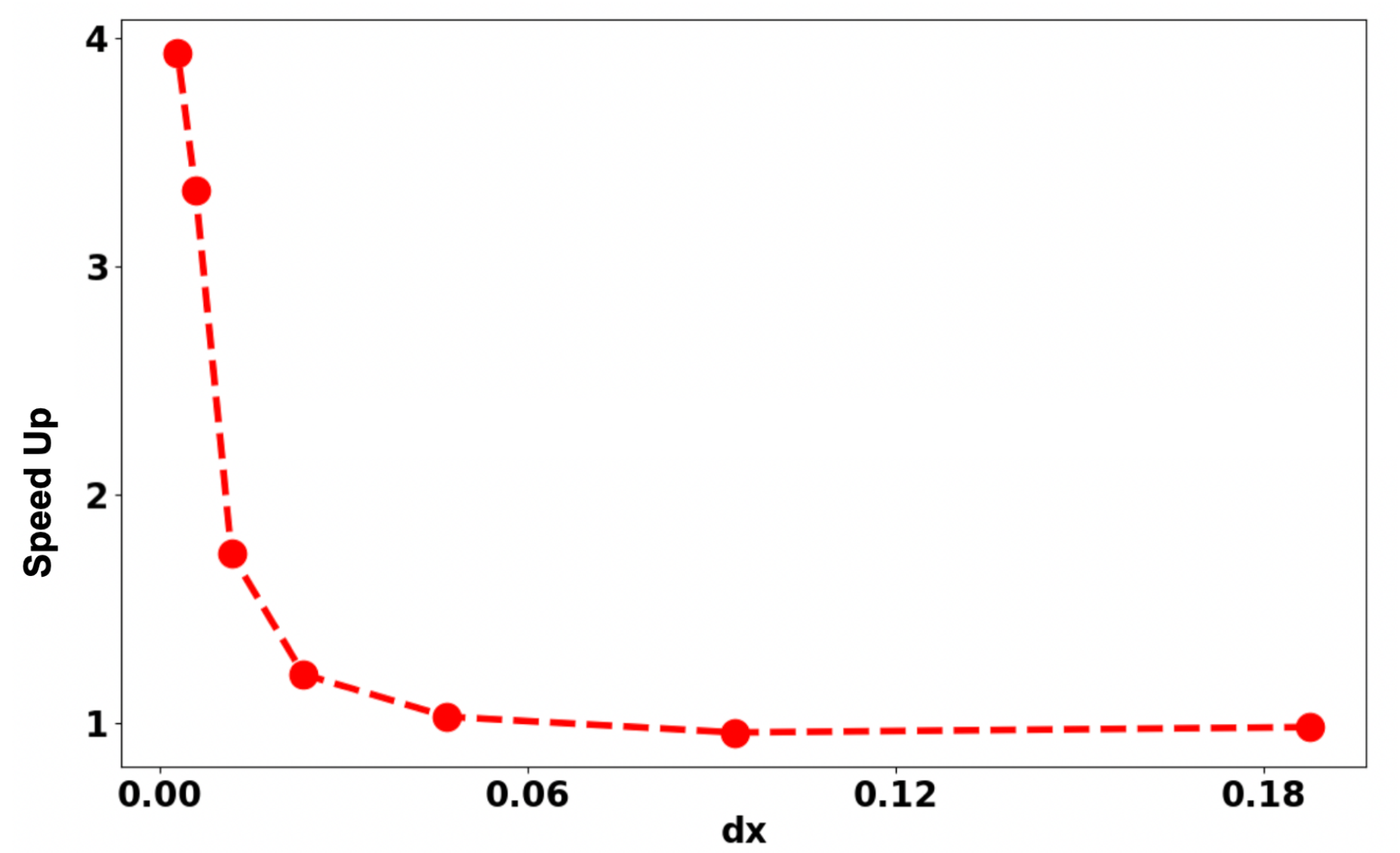}
	\caption{Comparison of execution times between the transformed non-uniform approach and the equivalent uniform grid solution.}
	\label{time_comparison}
\end{figure}

\vspace{3mm}
Many physical systems exhibit localized regions of high field variation or strong source gradients that require fine spatial resolution for accurate numerical solutions. Examples include charged particle beams, plasma sheaths, sharp interfaces in fluid flow, and local hot-spots in heat conduction. Using a uniform grid with sufficiently small spacing to resolve these features throughout the entire domain can be computationally expensive, leading to excessive memory usage and long execution times. To overcome this, we apply our transformation-based approach to design a non-uniform grid that is refined only where needed, while remaining coarse in regions where the solution varies slowly.\\
In this demonstration, we consider the two-dimensional Poisson equation with a Gaussian charge density localized at the center of the domain. The computational domain is initially defined as a uniform square grid of size $\mathrm{L_x = L_y = 6}$, with $\mathrm{N_x = N_y = 129}$ grid points. This uniform grid is then transformed into a non-uniform physical grid via the following coordinate mapping:
\begin{equation}
    \label{eq52}
    \begin{split}
        \mathrm{x^{\prime}} &= \mathrm{s \sinh(\alpha x)}, \\
        \mathrm{y^{\prime}} &= \mathrm{s \sinh(\alpha y)},
    \end{split}
\end{equation}
where the parameters $\mathrm{s = 5}$ and $\mathrm{\alpha = 0.6}$ control the overall stretching scale and the rate of refinement, respectively. The $\sinh$ transformation has the desirable property of producing higher grid density near the origin ($x, y \rightarrow 0$), while smoothly coarsening the mesh away from the center. Figures~\ref{stretch_grid}(a) and (b) show the original uniform computational grid and the resulting non-uniform physical grid.\\
On this stretched grid, we impose a Gaussian charge density beam defined as:
\begin{equation}
\label{eq53}
\mathrm{\rho' = -\rho_0 \exp\left( -\frac{(x^{\prime} - x_0)^2}{2\sigma_x^2} - \frac{(y^{\prime} - y_0)^2}{2\sigma_y^2} \right)},
\end{equation}
where $\mathrm{x_0 = y_0 = 0.0}$ places the beam center at the origin, $\mathrm{\rho_0 = 100}$ sets its peak amplitude, and $\mathrm{\sigma_x = \sigma_y = 0.5}$ control the beam width. The non-uniform grid naturally allocates more points in the beam region, improving resolution without increasing the total number of grid points excessively.\\
For comparison, Fig.~\ref{stretch_grid}(c) shows the equivalent uniform grid that would be required to achieve the same smallest grid spacing as the stretched grid. This resolution requirement is calculated as:
\[
\mathrm{N_{equiv}} \approx \frac{\mathrm{Domain \ size \ of \ physical \ stretched \ grid}}{\mathrm{Minimum \ grid \ spacing\  on \ physical \ stretched \ grid}},
\]
which in this case corresponds to a uniform grid of approximately $513 \times 513$ points — about 16 times more points than our transformed approach.\\
To solve the Poisson equation, we compute the Jacobian matrix of the transformation from Eq.~\ref{eq52}, which allows evaluation of the anisotropic permittivity tensor $\boldsymbol{\varepsilon}$ used in the transformed formulation. The charge density is mapped back to the uniform computational grid via $\mathrm{\rho = \rho^{\prime} \det(\mathcal{J})}$, where $\mathcal{J}$ is the Jacobian. In this particular transformation, the off-diagonal components vanish ($\varepsilon_{xy} = \varepsilon_{yx} = 0$), and the diagonal terms are shown in Fig.~\ref{epsilon_stretch}.\\
Since the potential remains invariant under the coordinate transformation, we can solve the transformed Poisson equation directly on the original $129 \times 129$ uniform computational grid and then map the resulting solution back to the physical grid. Figure~\ref{potential_stretch}(a) shows the potential distribution on the non-uniform grid, while Fig.~\ref{potential_stretch}(b) shows the equivalent solution obtained on the much larger $513 \times 513$ uniform grid.\\
Finally, the computational benefit of the approach is quantified in Fig.~\ref{time_comparison}, which compares execution times for increasing grid sizes. While the performance difference is small at lower resolutions, the advantage becomes substantial for larger grids, with our method achieving similar accuracy at approximately one quarter of the computational cost.
\subsection{Potential flow over an immersed body} 
\begin{figure}[H]
	\centering
	\includegraphics[width=1.0\textwidth]{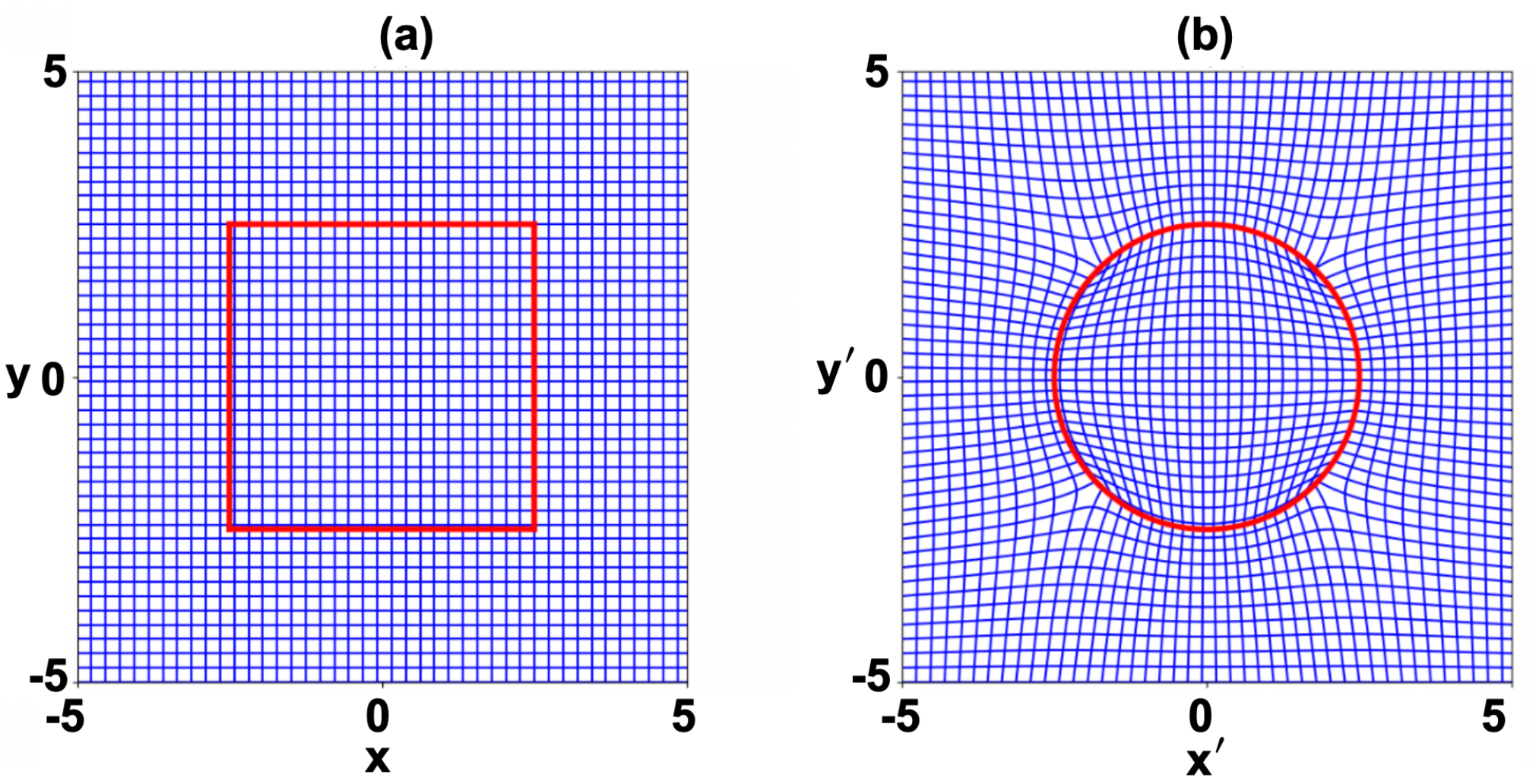}
	\caption{(a) Computational domain with selected transformation region (red outline), and (b) corresponding transformed physical domain.}
	\label{flow_1_grid}
\end{figure}

\begin{figure}[H]
	\centering
	\includegraphics[width=1.0\textwidth]{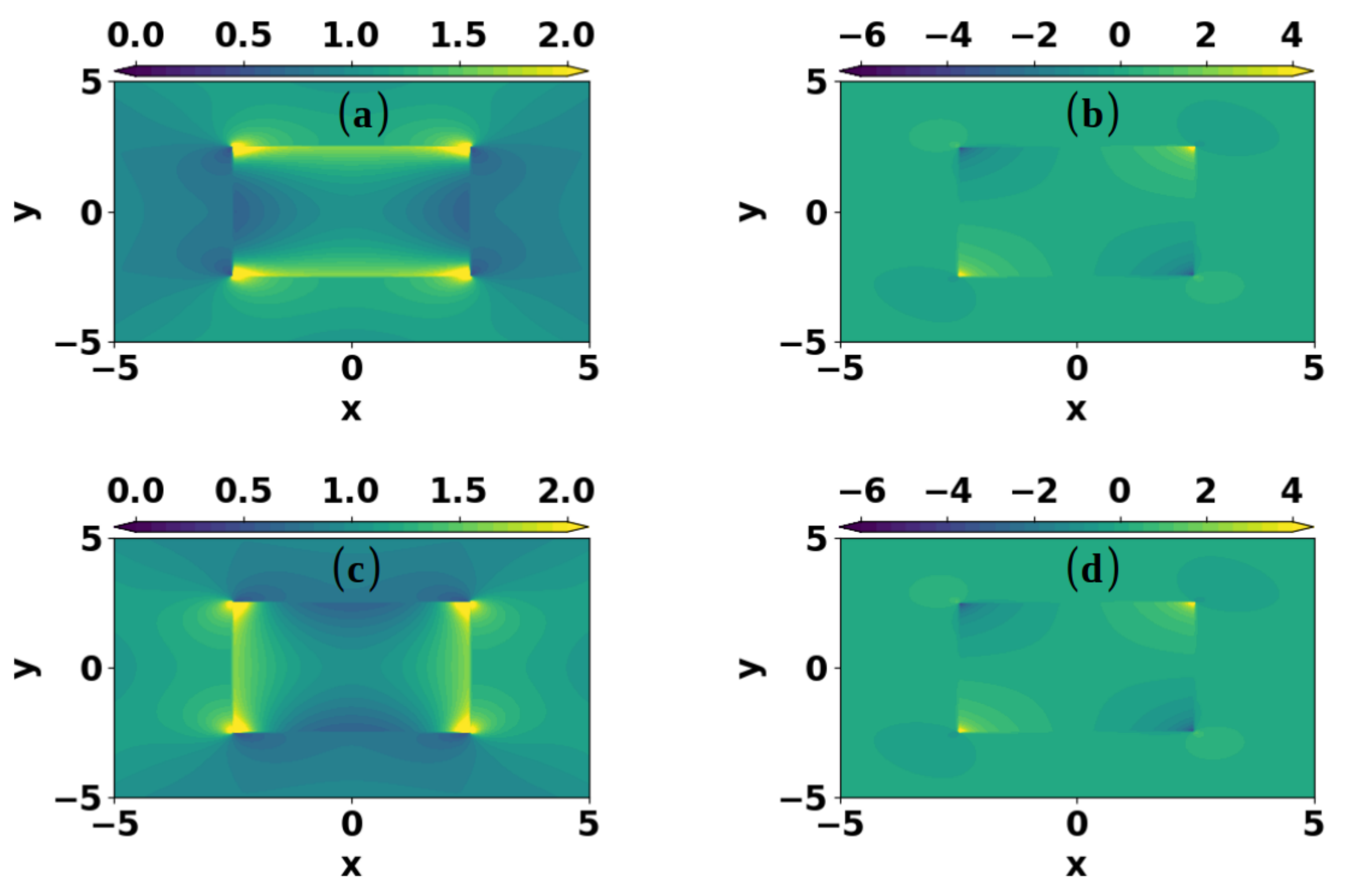}
	\caption{Components of the permittivity tensor: (a) $\varepsilon_{xx}$, (b) $\varepsilon_{xy}$, (c) $\varepsilon_{yy}$, and (d) $\varepsilon_{yx}$.}
	\label{flow_1_eps}
\end{figure}

\begin{figure}[H]
	\centering
	\includegraphics[width=1.0\textwidth]{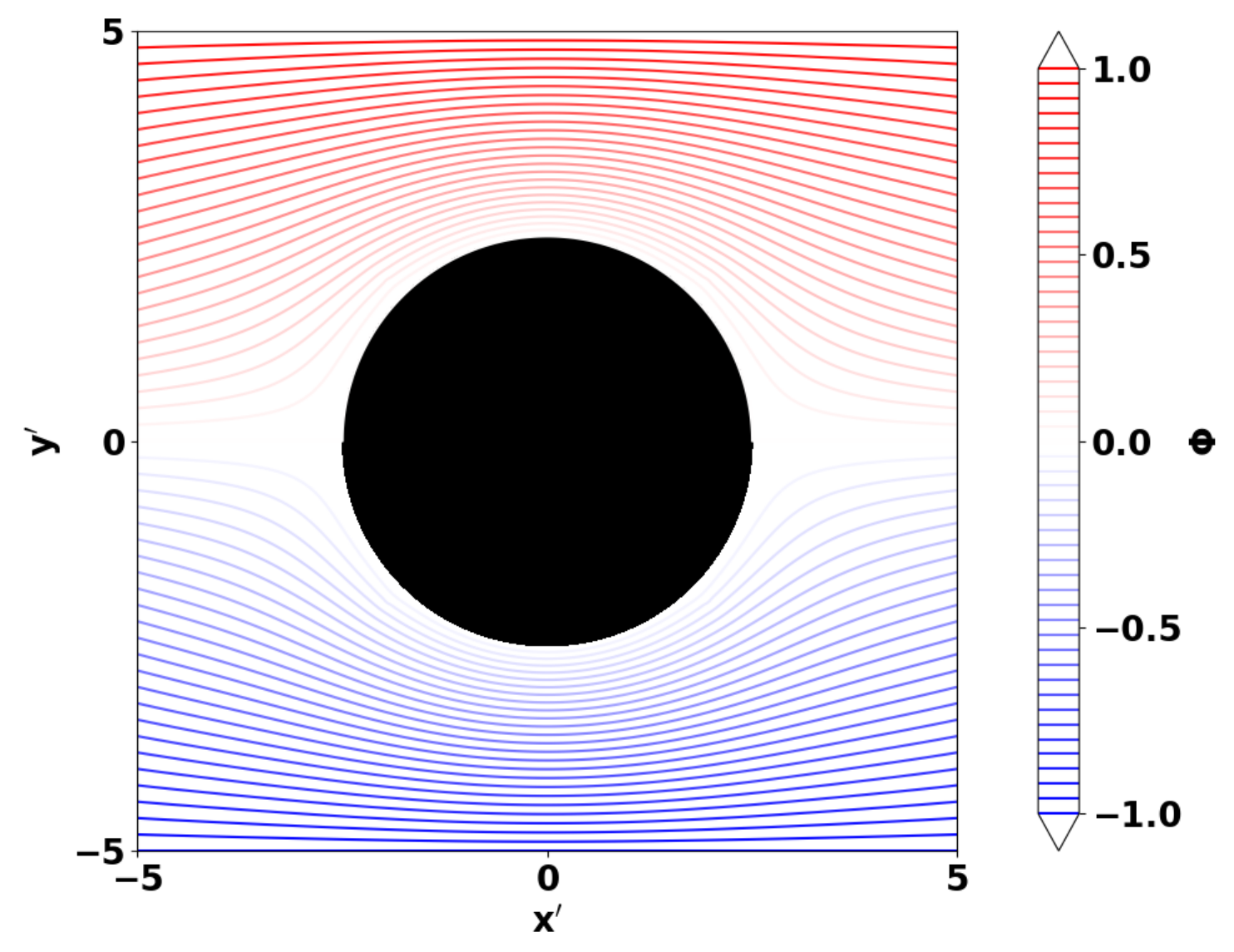}
	\caption{Potential flow around a stationary circular cylinder.}
	\label{flow_1_sol}
\end{figure}

\begin{figure}[H]
	\centering
	\includegraphics[width=1.0\textwidth]{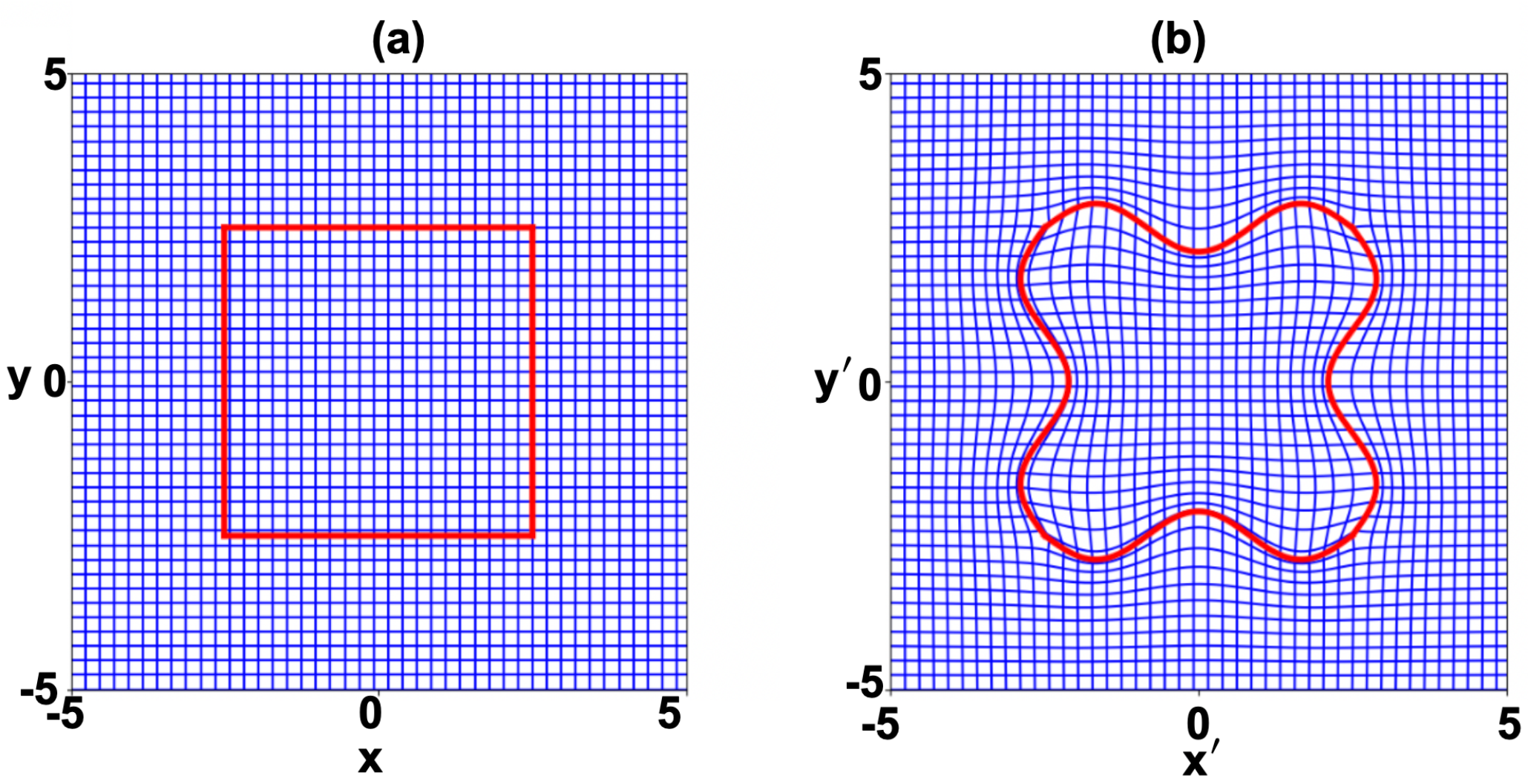}
	\caption{Computational and physical domains for the arbitrary-shape obstacle. The boundary deformation is confined to a central subregion, leaving the far-field mesh undistorted.}
	\label{flow_arb_grid}
\end{figure}

\begin{figure}[H]
	\centering
	\includegraphics[width=1.0\textwidth]{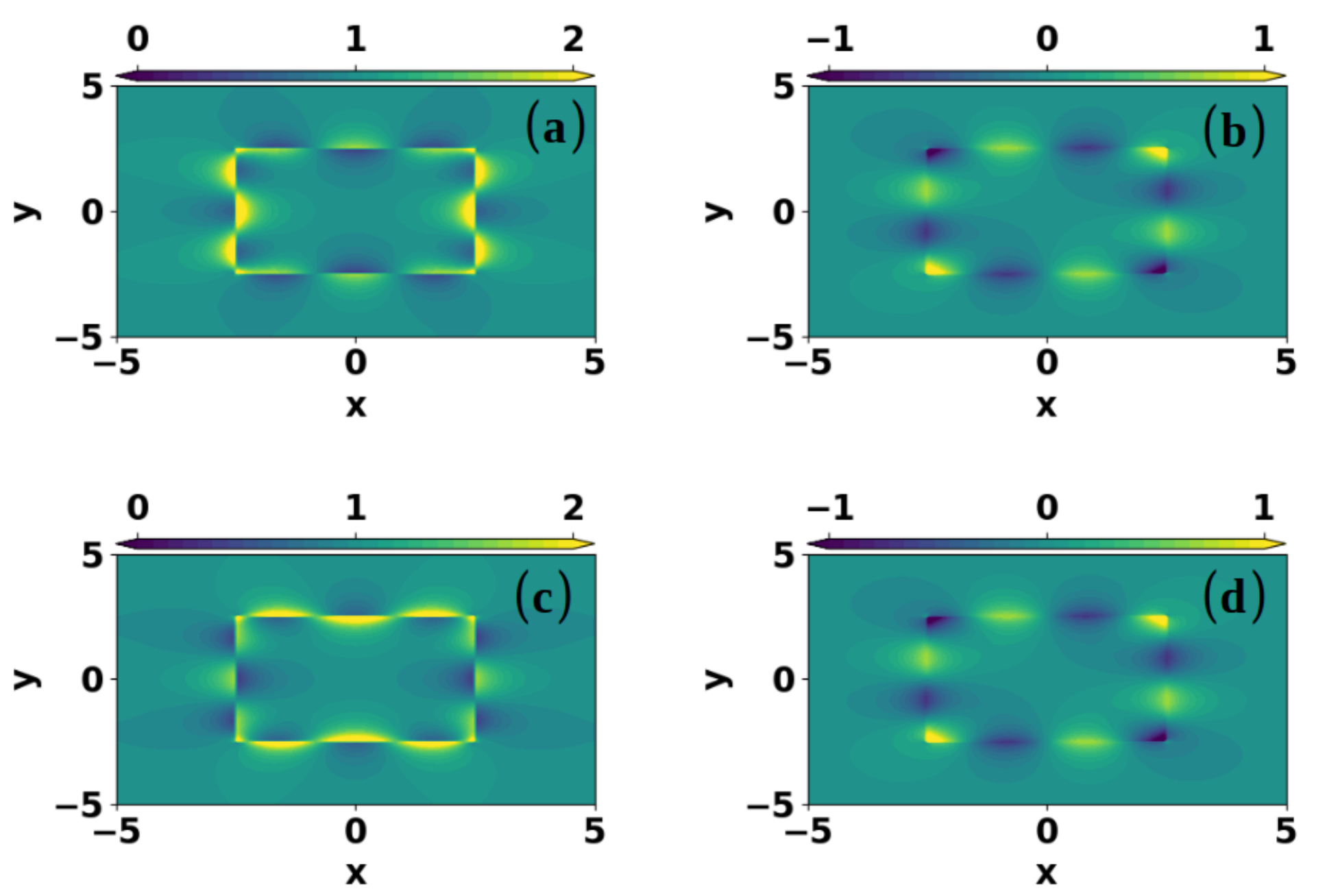}
	\caption{Components of the anisotropic permittivity tensor for the arbitrary-shape transformation. Strong variations occur near the perturbed boundary.}
	\label{flow_arb_eps}
\end{figure}

\begin{figure}[H]
	\centering
	\includegraphics[width=1.0\textwidth]{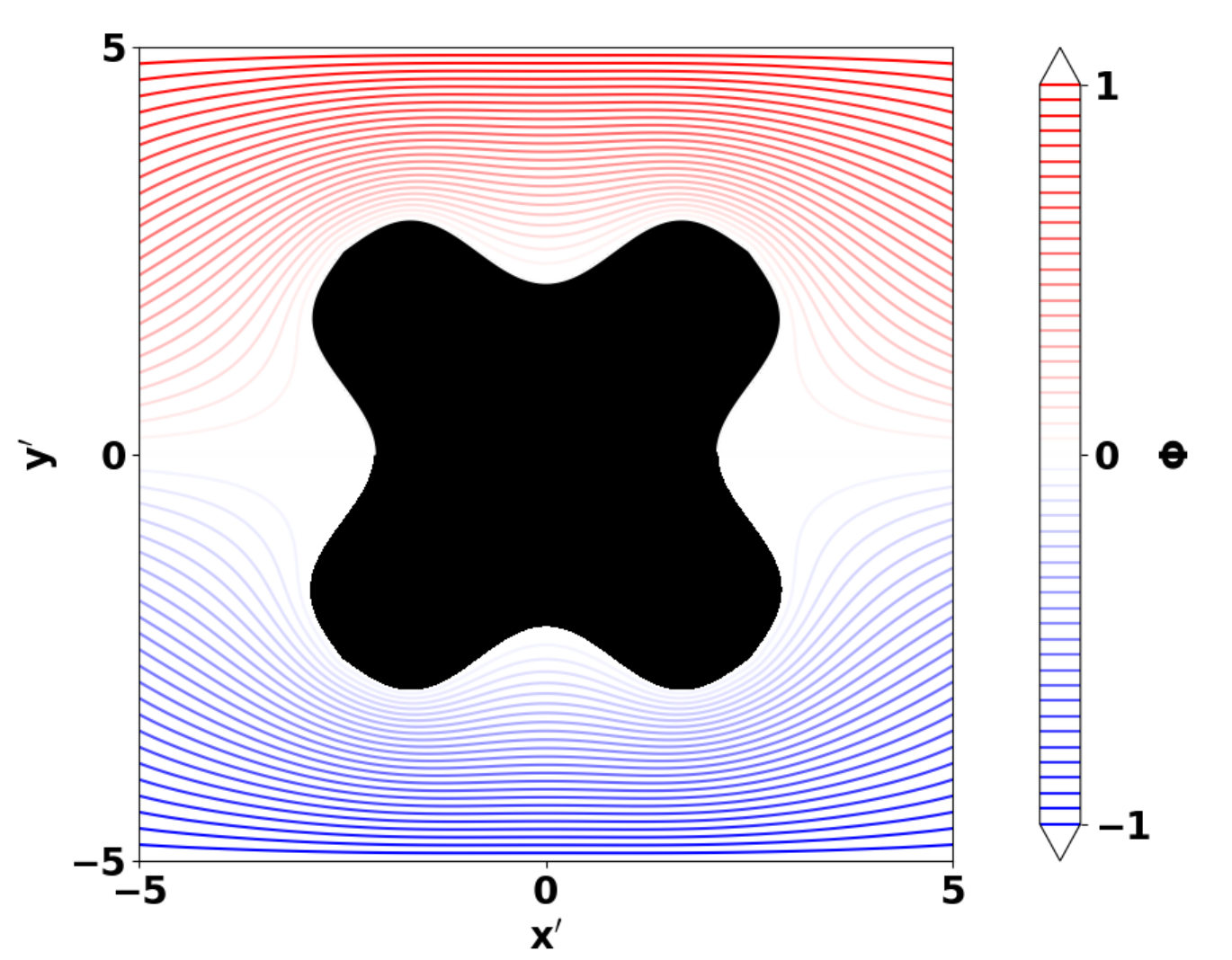}
	\caption{Simulated potential flow around the arbitrary-shaped obstacle.}
	\label{flow_arb_sol}
\end{figure}

\vspace{5mm}
 The flow of an ideal (incompressible, inviscid and irrotational) fluid around a circular cylinder is described by a stream function, which satisfies the Laplace equation. The level-surfaces of this function represent the streamlines of the flow field. In this subsection, we present the flow field contours for this problem using our method. For such problems, the boundary condition needs to be applied on an irregular surface placed inside the main domain of the simulation. Apart from fluid dynamics, a similar situation is encountered when one is interested in the electric field around the arbitrary shaped electrode. \\ As a first case for demonstration, we consider a classic problem of potential flow over a cylinder. We use identical mapping used in Eq.~\ref{eq44} to transform a square domain in the computational space to a circular domain in the physical space. It should be noted that in the present case, both square and circle are embedded inside the wider rectangular domains as can be seen from Fig.\ref{flow_1_grid}. In this figure, the boundaries of these embedded shapes are shown by the red color. Recall that the physical grid (Fig.\ref{flow_1_grid}(b)) is obtained by solving the Laplace equation over the computational grid (Fig.\ref{flow_1_grid}(a)) with the boundary condition described by the mapping Eq.~\ref{eq44}. Ideally, more sophisticated grid generation algorithms can also be used for this purpose. The components of the modified permittivity tensor used on the computational grid are in Fig.~\ref{flow_1_eps}. 
The boundary conditions correspond to a uniform incident/exiting flow in the $x-$ direction. Therefore, the stream function is made to vary linearly along the lines $x=0$ and $x=L_x$. Additionally, the stream function is held constant along the top ($y = L_y$) and bottom ($y = 0$) edges of the simulation domain, corresponding to undisturbed flow away from the cylinder. Finally, the stream function is also held constant along the surface of the cylinder, corresponding to the no-penetration boundary condition. Note that in the computation domain, the no-penetration boundary condition can be trivially applied on the embedded red-colored square, since it is exactly aligned with the grid. These boundary conditions are summarized as follows:
\begin{equation}
\label{eq47}
\begin{split}
\phi(0,y)     &= \phi_0 + \frac{\phi_1 - \phi_0}{L_y} \, y, \\
\phi(L_x,y)   &= \phi_0 + \frac{\phi_1 - \phi_0}{L_y} \, y, \\
\phi(x,0)     &= \phi_0, \\
\phi(x,L_y)   &= \phi_1, \\[3pt]
\phi(x,y)     &= 0, \quad \text{for} \ (x, y) \in \Omega_{\text{cylinder}},
\end{split}
\end{equation}
In this study, we set $\phi_0 = -1$ and $\phi_1 = 1$. The resulting streamlines are shown in Fig.~\ref{flow_1_sol}. These results are consistent with the theoretical results of the flow past a cylinder. The flow field can be constructed by taking the derivative of the streamfunction. Note that an identical mapping with modified boundary condition on $\Phi$ can be used to calculate the electric field lines emerging from the circular electrode. 
In many engineering and physical applications, obstacles do not conform to simple geometric primitives such as circles or ellipses. Instead, complex or irregular boundaries must be accommodated, often arising from natural formations, optimized aerodynamic profiles, or manufacturing constraints. To demonstrate the flexibility of the transformation-based approach, we consider here a two-dimensional potential flow around an arbitrary-shaped obstacle generated via controlled boundary deformation.\\
The arbitrary geometry is created by prescribing a sinusoidal perturbation along the boundaries of a selected rectangular subdomain $[x_0, x_c] \times [y_0, y_c] \subset \Omega$. The transformation for the physical coordinates $(x', y')$ is defined as:
\begin{equation}
\label{eq51}
\begin{aligned}
\mathrm{x^{\prime}} &= 
\begin{cases}
\mathrm{x - a \Delta x \,\sin\left( \pi n\dfrac{y - y_{0}}{y_{c} - y_{0}} \right)}, & \mathrm{x = x_0,\ y_0 \leq y \leq y_c}, \\
\mathrm{x + a \Delta x \,\sin\left( \pi n\dfrac{y - y_{0}}{y_{c} - y_{0}} \right)}, & \mathrm{x = x_c,\ y_0 \leq y \leq y_c}, \\
\mathrm{x}, & \text{otherwise},
\end{cases}
\\[1em]
\mathrm{y^{\prime}} &= 
\begin{cases}
\mathrm{y - a \Delta y \,\sin\left( \pi n\dfrac{x - x_{0}}{x_{c} - x_{0}} \right)}, & \mathrm{y = y_0,\ x_0 \leq x \leq x_c}, \\
\mathrm{y + a \Delta y \,\sin\left( \pi n\dfrac{x - x_{0}}{x_{c} - x_{0}} \right)}, & \mathrm{y = y_c,\ x_0 \leq x \leq x_c}, \\
\mathrm{y}, & \text{otherwise}.
\end{cases}
\end{aligned}
\end{equation}
In this formulation, $a$ controls the amplitude of the deformation, $n$ determines the number of oscillations along each perturbed edge, and $\Delta x$, $\Delta y$ represent the computational grid spacing. The sinusoidal form ensures smooth transitions between peaks and troughs, while keeping the boundary periodicity well-defined. This approach enables systematic control of boundary complexity without requiring an explicit parametric definition of the entire shape.\\
To generate the interior deformation corresponding to these boundary modifications, Laplace’s equations for $x'(x,y)$ and $y'(x,y)$ are solved across the computational domain with the transformed boundary conditions given in Eq.~\ref{eq51}. This harmonic interpolation produces a smooth, grid-conforming mapping from the computational domain to the physical domain, even for strongly deformed boundaries.\\
For the simulation, the domain size is set to $\mathrm{L_x = L_y = 10}$ with a uniform computational grid resolution of $\mathrm{N_x = N_y = 257}$. The transformation parameters are chosen as $\mathrm{a = 10}$ and $\mathrm{n = 3}$, with $\Delta x = \dfrac{L_x}{N_x - 1}$ and $\Delta y = \dfrac{L_y}{N_y - 1}$. The inner transformed domain is specified by $\mathrm{x_0 = -2.5}$, $\mathrm{x_c = 2.5}$, $\mathrm{y_0 = -2.5}$, and $\mathrm{y_c = 2.5}$. The computational and physical grids are shown in Fig.~\ref{flow_arb_grid}. The smooth deformation of the mesh inside the inner region accurately follows the prescribed boundary oscillations, while the outer mesh remains regular.\\
The Jacobian of the transformation is then used to compute the modified anisotropic permittivity tensor, shown in Fig.~\ref{flow_arb_eps}, which is required to solve the Poisson equation on the uniform computational grid. The components of the computed permittivity tensor clearly indicate regions of significant grid stretching and skewing near the deformed boundary. As in the previous test cases, the same potential flow boundary conditions described in Eq.~\ref{eq37} are applied. Finally, the computed potential flow around the arbitrary shape is presented in Fig.~\ref{flow_arb_sol}, illustrating the method’s ability to capture complex flow patterns around non-standard geometries without resorting to unstructured meshes or excessive local refinement.\\

\section{Conclusion}
\vspace{5mm}
In this work, we have presented a general and flexible method for solving the Poisson/Laplace equation on complex geometries using a finite difference framework. The approach achieves second-order accuracy in the $L_{1}$, $L_{2}$, and $L_{\infty}$ norms and is applicable to both analytically and numerically defined coordinate transformations. By incorporating the effects of the coordinate mapping into an equivalent anisotropic permittivity tensor, the problem can be solved efficiently on a uniform computational grid while accurately representing highly non-uniform or irregular physical domains. This formulation also provides the freedom to apply efficient solvers such as the Geometric Multigrid (GMG) method directly on the uniform computational grid.\\
Several applications have been demonstrated, including the resolution of localized charge density beams on stretched grids. These examples highlight that the proposed method can achieve the resolution of a much finer uniform grid—such as replacing a $513 \times 513$ uniform grid with a $129 \times 129$ computational grid—while significantly reducing computational cost. The technique is not limited to static geometries and can be applied to a wide range of complex domain shapes.\\
The method can be extended in several directions. First, it can be applied to three-dimensional problems, such as the 3D Poisson equation with Dirichlet boundary conditions, by computing the full $3 \times 3$ permittivity tensor arising from a 3D coordinate transformation. Second, higher-order accuracy can be achieved by incorporating higher-order finite difference schemes into the transformed formulation. Third, when the domain geometry evolves over time, the method remains applicable by updating the transformation and the associated permittivity tensor at each time step. Finally, by using an appropriately defined source term in the Laplace equation to generate a numerical coordinate transformation, it is possible to design highly complex and adaptive grids tailored to specific problem features.\\
In future work, we will be incorporating this algorithm in a two-dimensional electrostatic particle-in-cell (PIC) code. This will open-up possibilities of performing kinetic plasma simulations with complex geometries using a simple finite-difference code. Currently, more sophisticated codes, employing either finite element or finite volume methods, are required for such problem. 
             
\section{Acknowledgments}
\vspace{5mm}
Research work presented in this paper is fully supported by funding from UM-DAE Centre for Excellence in Basic Sciences (CEBS).


\bibliographystyle{iopart-num}
\bibliography{references}

\end{document}